\documentclass{amsart}
\usepackage{amscd,amssymb}
\usepackage[dvips]{graphicx}

\numberwithin{equation}{section}

\newtheorem{thm}{Theorem}
\newtheorem{co}[thm]{Corollary}

\newtheorem{pro}[thm]{Proposition}
\theoremstyle{remark}

\newtheorem{rem}[thm]{Remark}
\newtheorem{defi}[thm]{Definition}
\newtheorem{exa}[thm]{Example}
\newtheorem{nota}[thm]{Notation}

\begin{document}
\title{Coding and tiling of Julia sets 
for subhyperbolic rational maps }
\author{Atsushi KAMEYAMA}
\address{Division of Mathematical Science for Social Systems\\
Department of Systems Innovation\\
Graduate School of Engineering Science\\
Osaka University, Toyonaka, Osaka, 560-8531, Japan}
\email{kameyama@sigmath.es.osaka-u.ac.jp}
\urladdr{http://www.sigmath.es.osaka-u.ac.jp/\~{}kameyama/}
\date{}

\maketitle

\begin{abstract}
 Let $f:\hat{\mathbb C}\to\hat{\mathbb C}$ be a subhyperbolic  
rational map of degree $d$.
We construct a set of coding maps $\mathrm{Cod}(f)=
\{\pi_r:\Sigma\to J\}_r$ of 
the Julia set $J$ by 
geometric coding trees, where the parameter $r$ ranges over 
mappings from a certain tree to the Riemann sphere.
Using the universal covering space $\phi:\tilde S\to S$ 
for the corresponding orbifold, 
we lift the inverse of $f$  to an iterated function system 
$\mathcal I=(g_i)_{i=1,2,\dots,d}$. 
For the purpose of studying the structure of $\mathrm{Cod}(f)$, 
we generalize Kenyon and Lagarias-Wang's results 
: If the attractor $K$ of $\mathcal I$ 
has positive measure, then $K$ tiles $\phi^{-1}(J)$, and the 
multiplicity of $\pi_r$ is well-defined.
Moreover, we see that the equivalence relation induced by $\pi_r$ 
is described  by a finite directed graph, and give a necessary and 
sufficient condition for two coding maps $\pi_r$ and $\pi_{r'}$ to 
be equal.
\end{abstract}

\section{Introduction}

The method of symbolic dynamics is prevalent in the study
 of dynamical systems.
Particularly, to investigate an attractor (or repeller) of 
a dynamical system,
one often uses symbolic dynamics to code the attractor.
In the present paper, we study coding maps of Julia sets for rational 
maps.
We emphasize that we treat not only individual coding maps but 
totalities of coding maps.

We say a rational map $f:\hat{\mathbb{C}}\to\hat{\mathbb{C}}$ 
of the Riemann sphere to itself is subhyperbolic if each critical point 
is either preperiodic or attracted to an attracting cycle.
We will construct a set of coding maps $\mathrm{Cod}(f)=\{\pi_r\}_r$ 
from the full 
shift to the Julia set $J$ by using `geometric coding trees.' 
Geometric coding tree technique is developed by Przytycki and 
his coauthors for 
general holomorphic maps (\cite{Prz85}, \cite{Prz86}, \cite{PrSk91}).
See also \cite{Ly86}, \S 1.16.
The parameter $r$ ranges over  {\em radials}, 
which are mappings from 
a certain topological tree to the Riemann sphere. 
One goal of considering $\mathrm{Cod}(f)$ is to understand the
combinatorics of $f$.
In view of the construction of geometric coding trees, the structure 
of $\mathrm{Cod}(f)$ reflects the combinatorics of $f$.

We know that coding is effectively used in the study of 
dynamics of interval maps (see for example \cite{CoEc80}, \cite{MiTh88}, 
\cite{MeSt93}).
Coding often works in the parameter space as well as in the
dynamical space.
Recall that the dynamical space is just an interval $I$ with natural 
partition, that is, $I$ is divided 
into several subintervals by the turning points.
This partition gives a natural coding of $I$, and we obtain a nice 
invariant called the kneading sequence, which is 
defined as the symbol sequences corresponding to the forward orbits 
of the turning points.
For example, in a certain family of real polynomial maps the 
kneading sequences almost completely classifies these maps 
up to topological conjugacy. 
In particular, for the quadratic family $x\mapsto ax(x-1)$ 
($0<a<4$), we have the monotonicity of the kneading sequence (and 
the topological entropy).   
Roughly speaking, the natural coding parametrizes the bifurcation of the 
quadratic family. 
However, in a larger family of complex rational maps, 
coding does not seem to work well in the parameter space.
The main reason of this difficulty is the absence of natural partitions.
We do not have a nice invariant like the kneading sequence.
In such a situation, thus it is less important to consider individual 
coding maps.
This is why we treat totalities of coding maps.

A complete description of $\mathrm{Cod}(f)$ is quite difficult except 
a few cases including $f(z)=z^2$, $f(z)=z^2-2$, etc.
It is unfortunate that even the case $f(z)=z^d$ with $d\ge 3$ 
and the Cantor set case (e.g. $f(z)=z^2-3$) are complicated.
Thus we will try to find  tools to manage
$\mathrm{Cod}(f)$, keeping in mind the following natural and 
naive problems: (1) What is the canonical coding?
(2) Are there any good structures on $\mathrm{Cod}(f)$?
The present paper does not completely solve these problem, 
but gives several fundamental facts which will be a help to approach 
these problems.
Our main results are concerned with the 
multiplicities of coding maps and the equivalence relations on the 
space of symbol sequences.

In our setting, $f^{-1}$ can be lifted by the `universal 
covering' $\phi:\tilde S\to S=\hat{\mathbb C}-\mathrm{AP}$, 
where $\mathrm{AP}$ is the set of attracting periodic points.
Let $d$ be the degree of $f$.
For $i=1,2,\dots,d$, there exists a holomorphic contraction
$g_i:\tilde S\to \tilde S$ depending on $r$ such that the diagram

$$\begin{CD}
   \tilde S-\phi^{-1}f^{-1}(\mathrm{AP}) @<{g_i}<< \tilde S\\
   @V\phi VV         @VV\phi V \\
   S-f^{-1}(\mathrm{AP})       @>> f >       S
  \end{CD}
$$
commutes.
Thus there exists a compact set $K$ such that $K=\bigcup_{i=1}^d 
g_i(K)$ by Hutchinson \cite{Hu81}.
We call $K$ the {\em Julia tile} with respect to $r$.
If $r$ is `suitable', then the coding map $\pi_r$ is onto, and so 
$\phi(K)=J$. 

The situation above is analogous to that of self-affine tiling.
We recall self-affine tiling briefly.
Let $A$ be an $n\times n$ expanding integral matrix, and let 
$d=|\det A|$.
Since $A\mathbb{Z}^n\subset\mathbb{Z}^n$, its projection $f:\mathbb T^n
\to\mathbb T^n$ is well-defined, where $\mathbb T^n=\mathbb R^n/\mathbb 
Z^n$ is the $n$-torus.
A lift of $f$ has the form $x\mapsto Ax+a,a\in \mathbb Z^n$.
Choosing $a_1,a_2,\dots,a_d\in\mathbb Z^n$, we obtain $d$ contractions
$g_i:x\mapsto A^{-1}(x-a_i)$ by which we have the commuting diagram 
$$\begin{CD}
   \mathbb R^n @<{g_i}<< \mathbb R^n\\
   @V\phi VV         @VV\phi V \\
   \mathbb T^n       @>> A >       \mathbb T^n
  \end{CD}
$$
for $i=1,2,\dots,d$, and the self-affine set $K$ with $K=\bigcup_{i=1}^d 
g_i(K)$.
Kenyon \cite{Ken92} showed that if $K$ has positive Lebesgue measure, 
then $K$ tiles $\mathbb R^n$ (i.e. there exists a set of translations 
$T\subset \mathbb Z^n$ such that (1) $\mathbb R^n=\bigcup_{t\in T}(K+t)$, 
and (2) the Lebesgue measure of $(K+t)\cap (K+t')$ vanishes for distinct 
 $t,t'\in T$).
See for more details \cite{Ken92} and \cite{LaWa96b}.

We generalize Kenyon and Lagarias-Wang's results and obtain:

\medskip

\noindent{\bf Tiling Theorem} (Theorem \ref{thm:2}).
Let $\mu$ be the equilibrium state for $f:J\to J$ with constant 
weight, and $\tilde\mu$ the lift of $\mu$ to $\tilde S$.
If $\tilde\mu(K)$ is positive, then $K$ tiles $\tilde J=\phi^{-1}(J)$, 
that is, there exists a set $T$ of deck 
transformations of $\phi:\tilde S\to S$ such that  
 (1) $\tilde J=\bigcup_{t\in T}t(K)$, 
and (2) $\tilde\mu(t(K)\cap t'(K))=0$ for distinct  $t,t'\in T$.

\medskip

\noindent{\bf Multiplicity Theorem} (Theorem \ref{thm:21}).
If $\tilde\mu(K)$ is positive, then there exists $n>0$ depending on 
$r$ such that $\pi_r$ is {\em almost $n$-to-one}, 
that is, $\#\pi_r^{-1}\pi_r(\omega)=n$ for $m$-almost 
all $\omega\in\Sigma$, where $m$ is the Bernoulli measure with
identically distributed weight. 
The number $n=n_r$ is said to be the multiplicity of $\pi_r$.
\medskip

For a given coding map $\pi_r$, it is difficult to calculate its 
multiplicity in general.
The most primitive way to do that is directly seeing the equivalence 
relation induced by $\pi_r$ on $\Sigma$.
We will show that the equivalence relation is described by a finite 
graph, which is a version of Fried's result (\cite{Fri87}, Lemma 1).
See also \cite{Ka03}.

\medskip 

\noindent{\bf Finite Graph Theorem} (Theorem \ref{thm:28}). 
We can construct a directed graph with vertex set $V$, edge set $E$, and 
weight $\alpha:E\to \{1,2,\dots,d\}^2$ such that 
$\pi_r(\omega_1\omega_2\cdots)=\pi_r(\omega'_1\omega'_2\cdots)$ 
if and only if there exists a sequence $e_1,e_2,\ldots\in E$ with 
$\alpha(e_i)=(\omega_i,\omega'_i)$ and $(\text{the terminal vertex of 
$e_i$})=(\text{the initial vertex of $e_{i+1}$})$ for $i=1,2,\dots$.

\medskip

We expect that any subhyperbolic rational maps have the canonical 
coding maps.
The word `canonical' has vagueness, but it become clearer by the notion 
of multiplicity.
If there exists the simplest nontrivial coding map, we may consider it 
canonical.
We will look at several examples of subhyperbolic rational maps later,
and will find out that each of them has a nontrivial 
coding map which is apparently the simplest.
These coding maps have the features (1) the multiplicities are equal to 
one, (2) the {\em connecting sets} 
$$E=\bigcup_{i\ne j}(g_i(K)\cap g_j(K))
\cup\bigcup_i \{x\in g_i(K)\,|\,\#g_i^{-1}(x)\ge 2\}$$
are small, and (3) the Julia tiles $K$ have simple shapes.
While the second and the third features are still vague, the first one 
mathematically makes sense.
Thus our conjecture is that any subhyperbolic rational maps have a  
coding maps with multiplicity one.

Another conjecture is that $\mathrm{Cod}(f)$ has some structure.
For example, we expect that there exists a natural action on 
$\mathrm{Cod}(f)$ by which we can control the diversity of multiplicities 
of $\pi_r$.

\medskip
\noindent{\bf Structure Theorem} (Theorem \ref{thm:24}).
$\pi_r=\pi_{r'}$ if and only if $r$ and $r'$ are freely homotopic mod 
$N_r$.

\medskip
This theorem says that 
$$\mathrm{Cod}(f)\approx
\bigcup_N\{r\,|\,N_r=N\}/(\text{freely homotopic mod $N$}),$$ 
where $N_r$ denotes the maximal invariant subgroup with respect to $r$ 
of the fundamental group of $\mathbb C-\{\text{postcritical points}\}$ 
we will define later.
Let $\mathcal A(f)$ be the monoid of rational maps commutative to $f$.
Then $\mathcal A(f)$ acts on $\mathrm{Cod}(f)$ by 
$(R,\pi)\mapsto R\circ\pi$, and the multiplicity of $R\circ\pi$ is 
equal to the degree of $R$ times the multiplicity of $\pi$ 
(Proposition \ref{pro:32}).

The organization of the present paper is as follows.
After giving the definitions of coding maps and Julia tiles in Section 
\ref{sec:def}, we consider several examples of subhyperbolic rational 
maps in Section \ref{sec:3}.
In one example, the L\'evy Dragon appears as a Julia tile.
Section \ref{sec:4} supplies the definitions of invariant subgroups 
and the equilibrium states, and shows some basic facts.
We prove Tiling Theorem and Multiplicity Theorem in Section 
\ref{sec:tile}.
In Section \ref{sec:st}, we discuss the structure of $\mathrm{Cod}(f)$, 
prove Structure Theorem, and show a couple of examples. 
Finite Graph Theorem  is proved in Section \ref{sec:5}.
In Section \ref{sec:nonhyp}, we see that several results hold for 
non-subhyperbolic rational maps.

\section{Definition}\label{sec:def}

In this section, we give the construction of coding maps after 
Przytycki.  
We obtain Julia tiles by lifting the coding maps to the universal 
covering spaces.
Let $f:\hat{\mathbb C}\to\hat{\mathbb C}$ be a rational map of degree 
$d$. 

\begin{defi}
We say that $$C=C_f=\{\text{critical points of $f$}\}
=\{c\,|\,\text{$f$ is not locally homeomorphic at $c$}
\}$$ is the {\em critical set} of $f$. 
The {\em postcritical set} is defined by 
$$P=P_f=\overline{\{f^k(c)\,|\,c\in C,n>0\}}.$$
Let $\mathrm{AP}$ denote the set of attracting periodic points of $f$.
We write $$S'=\hat{\mathbb{C}}-P\quad\text{and}\quad
 S=\hat{\mathbb{C}}-\mathrm{AP}.$$
\end{defi}

\begin{defi}
 Let $$Q=Q_d=\bigcup_{i=1}^d[0,1]_i/(0_i\sim 0_j)$$
be the topological tree made of $d$ copies of the unit interval 
with all the origins identified.
We say that a continuous map $r:Q\to S'$ is a 
{\em radial} if $r(0)=f\circ r(1_i)$ for $i=1,2,\dots,d$.
We say that $\bar x=r(0)$ is the {\em basepoint} of $r$.
A radial is considered as a $d$-tuple of curves $(l_i:[0,1]\to S')_i$ 
with the same initial point $\bar x$ such that $f\circ l_i(1)=\bar x$.
A radial $r$ is said to be {\em proper} if $r(1_i)\ne r(1_j)$ whenever
$i\ne j$.
We write the set of radials with basepoint $\bar x$ by 
$$\mathrm{Rad}(f,\bar x). $$
 Set 
$$\begin{array}{l}
L=L(f,{\bar x})=\{l:[0,1]\to S'\,|\,\text{$l$ is continuous and 
$l(0)=\bar x$}\},\\
\Lambda(f,\bar x)=\{l\in L(f,\bar x)\,|\,l(1)\in f^{-1}(\bar x)\}.
\end{array}$$
Then $\mathrm{Rad}(f,\bar x)=\Lambda(f,\bar x)^d$.

For a curve $l\in L$ and $x\in f^{-k}(\bar x)$, 
we define a curve $F_x(l):[0,1]\to S$ as the lift of $l$ by $f^k$ 
with initial point $x$, that is, $f^k\circ F_x(l)=l$ and $F_x(l)(0)=x$.  
Since $l$ does not pass through $P$, the curve $F_x(l)$ is uniquely 
defined.
\end{defi}

\begin{nota}
Let $(\sigma,\Sigma)$ be the one-sided fullshift of $d$ symbols.
Namely, $$\Sigma=\{1,2,\dots,d\}^{\mathbb{N}}$$ is the set of 
one-sided infinite sequences of $\{1,2,\dots,d\}$, and $\sigma:\Sigma\to
\Sigma$ is the shift map ($\sigma(w_1w_2\cdots)=w_2\cdots$).

Let $W$ be the set of words (i.e. finite sequences) of $d$ symbols, 
and $W_k$ the set of words of length $k$:
$$W=\bigcup_{k=1}^\infty W_k,\quad W_k=\{1,2,\dots,d\}^k.$$ 
For $w\in W$, we write 
$$\Sigma(w)=\{w\omega_1\omega_2\cdots\,|\,\omega_1\omega_2
\cdots\in\Sigma\}\subset\Sigma. $$
\end{nota}

Construction of coding maps.
Suppose $f$ is subhyperbolic.
For a radial $r=(l_i)_i$, we inductively define curves $l_w:[0,1]\to 
S'$ and points $x_w$ for $w\in W$.
First we set $x_i=l_i(1)$ for $i\in W_1=\{1,2,\dots,d\}$.
If $l_w$ and $x_w$ are determined for $w\in W_k$, 
we set $l_{iw}=l_i\cdot F_{x_i}(l_w)$ and 
$x_{iw}=l_{iw}(1)$ for $i\in W_1$.
By the expandingness of $f$, 
we have $l_\omega=\text{``$\lim_{k\to\infty}$''}
l_{\omega_1\omega_2\cdots\omega_k}$ and $x_\omega=\lim_{k\to\infty}
x_{\omega_1\omega_2\cdots \omega_k}$ for $\omega=\omega_1\omega_2\cdots
\in\Sigma$.
Note that $l_{\omega_1\omega_2\dots\omega_k}$ does not converge, but 
by a suitable change of parametrization, the limit $l_\omega$ is 
well-defined and unique up to parametrization.
Clearly, $x_\omega=l_\omega(1)$ for $\omega\in\Sigma$.
Since $x_\omega$ is an accumulation point of $f^{-k}(\bar
x),k=1,2,\dots$, the point belongs to the Julia set (for example, 
see \cite{Mi99}, 4.7).
It is easily seen that $f(x_\omega)=x_{\sigma\omega}$ and 
the mapping $\omega\mapsto x_\omega$ is continuous.
We denote this mapping by $\pi=\pi_r:\Sigma\to J$, and call it 
the {\em coding map} of $J$ for $r$.
We write $$\mathrm{Cod}(f)=\{\pi_r\,|\, \text{$r$ is a radial}\}.$$

\begin{rem}\label{rem:5}
\ 
\begin{itemize}
\item
$l_{u}\cdot F_{x_u}(l_w)$ is equal to $l_{uw}$ 
up to parametrization for $u\in W$ and $w\in W$ (or $w\in\Sigma$). 
\item $f\pi(\Sigma)=\pi(\Sigma)$.
\item If $r$ is proper, then $\pi_r:\Sigma\to J$ is onto.
However, the converse is not always true.

\item
 The image of $\pi$ is either a perfect set or a singleton.
Indeed, suppose $\pi(\Sigma)$ has a isolated point $p$.
Since $\pi^{-1}(p)$ is open, there exists $k$ such that 
$\sigma^k\pi^{-1}(p)=\Sigma$.
Therefore $\pi(\Sigma)=f^k(p)$.

\item
If $\pi(\Sigma)$ is perfect, then $\#f^{-k}(x)\cap\pi(\Sigma)\to\infty$ 
as $k\to\infty$ for any $x\in\pi(\Sigma)$.
Indeed, it is sufficient to show that there 
exists $k$ such that $\#f^{-k}(x)\cap\pi(\Sigma)\ge 2$.
Assume that $f^{-k}(x)\cap\pi(\Sigma)=\{y_k\}$ for any $k$.
There exists $k_0$ such that $y_k\notin C$ for $k\ge k_0$.
Take $\omega\in\Sigma$ with $\pi(\omega)=y_{k_0}$.
Then $\pi(w\omega)=y_{k_0+m}$ for any $w\in W_m$.
This means that $F_{x_w}(l_\omega)$'s are the same for $w\in W_m$.
Thus $x_w$'s are the same for $w\in W_m$.
Consequently, the accumulation points of $\{x_w\}_{w\in W_m}$ ($m\to
\infty$) consist of one point.

\end{itemize}
\end{rem}

From now on, we suppose that $f$ is subhyperbolic.

\begin{defi}\label{def:3}
A function $\rho:S\to\mathbb{N}$ is called a {\em ramification 
function} for $f$ 
if $\rho(x)=1$ for $x\notin P$ and $\rho(f(x))$ is a multiple 
of $\deg_xf\cdot\rho(x)$.
The minimal ramification function is called the {\em canonical} 
 ramification function and denoted by $\rho_f$. 

For a ramification function $\rho$, we have a {\em universal covering} 
$\phi:\tilde S\to S$ for the {\em orbifold} $(S,\rho)$ (i.e. 
$\tilde S$ is a connected and simply connected Riemann surface and 
$\phi$ is the holomorphic branched covering such that the local degree 
$\deg_{\tilde x} \phi$ is $\rho(\phi(\tilde x))$ for every 
$\tilde x\in S$).
See \cite{Mi99}, Appendix E.

 The universal covering $\phi:\tilde S\to S$ is constructed as 
follows.
Let $$G=\pi_1(S',\bar x)$$
be the fundamental group of $S'$.
For $$\gamma\in\Gamma=\Gamma(f,\bar x)
=\{\gamma\in L(f,\bar x)\,|\,\gamma(0)=\gamma(1)\},$$
we denote the homotopy class for $\gamma$ by $[\gamma]\in G$.
Let $\mathrm{AP}=\{a_1,a_2,\dots,a_p\}$ and $P-\mathrm{AP}=\{b_1,b_2,
\dots\}$.
Choose simple closed curves 
$A_1,A_2,\dots,A_{p},B_1, B_2,\dots\in\Gamma$ so that 
$A_i$ separates $a_i$ from the other points of $\mathrm{AP}$, and 
$B_j$ separates $b_j$ from the other points of $P$.
Then $G$ is generated by $[A_i],[B_j]\ (i=1,2,\dots,p,\ j=1,2,\dots)$.
Let $N=N^\rho$ be the normal subgroup of $G$ generated by 
$[B_j^{\rho(b_j)}],b_j\in P-\mathrm{AP}$. 
Then $[A_i^k]\notin N,k\in\mathbb Z$ and $[B_j^k]\notin N,
0<|k|<\rho(b_j)$ in our case.
The quotient group 
$$G^\rho=G/N$$
is called the fundamental group of the orbifold $(S,\rho)$.
 Let $$\bar L
=\{l:[0,1]\to S\,|\,\text{$l$ is continuous, $l(0)=\bar x$ 
 and $l(t)\in S'$
for $0\leq t<1$}\}.$$
We set $$\tilde S=\bar L/\sim_N,$$
where $l\sim_N l'$ if $l(1)=l'(1)$ and there exists a closed curve 
$\gamma\in \Gamma$ obtained by perturbing $ll'^{-1}$ near $l(1)$  
such that $[\gamma]\in N$. 
We obtain a branched covering $\phi:\tilde S\to S$ defined by 
$\phi([l])=l(1)$.
The surface $\tilde S$ is considered as either $\{|z|<1\}$ or 
$\mathbb{C}$.

The covering $\phi:\tilde S-\phi^{-1}(P)\to S'$ corresponds to 
the normal subgroup $N\subset G$.
Take $\tilde x\in\phi^{-1}(\bar x)$.
Since $(f\circ \phi)_*\pi_1(\tilde S-\phi^{-1}f^{-1}(P),z)
\supset N$, 
for $z\in \phi^{-1}f^{-1}(\bar x)$ 
there exists a holomorphic covering $g:\tilde S-\phi^{-1}(P)\to 
\tilde S-\phi^{-1}f^{-1}(P)$ such that 
\begin{equation}
f\circ \phi\circ g=\phi\label{eq:0}
\end{equation}
 and $g(\tilde x)=z$.
We extend $g$ to a holomorphic branched covering $g:\tilde S\to\tilde S$.
We say that $g$ is the {\em contraction associated with} $z$ (with 
respect to the basepoint $\tilde x$). 

By retaking a bigger ramification function if necessary, we may 
assume $\tilde S=\{|z|<1\}$. 
Thus $g$ is contracting with respect to the Poincar\'e metric.
Hence $f$ is uniformly expanding near $J$ with respect to the 
projection metric, that is, 
there exists $c>1$ such that $\|Df_z\|>c$ for $z$ in some neighborhood 
of $J$.
See \cite{Mi99}, \S19.
\end{defi}

\begin{defi}\label{def:8}
Fix $\tilde x\in \phi^{-1}(\bar x)$.
For $l\in \bar L$, we denote by $\tilde l$ 
the lift of $l$ to $\tilde S$ with $\tilde l(0)=\tilde x$
Let $r=(l_i)$ be a radial.
We have the contractions $g_1,g_2,\dots,g_d:\tilde S\to\tilde S$ 
associated with $\tilde l_i(1)$. 
Since $g_i$'s are contracting, we have the attractor $K$ of the iterated 
function system $\mathcal I=(g_1,g_2,\dots,g_d)$
(i.e. $K\subset \tilde S$ is the unique nonempty compact set 
with $K=\bigcup_{i=1}^dg_i(K)$. See \cite{Hu81}).
We use the notation 
$$g_w=g_{w_1}\circ g_{w_2}\circ\cdots\circ g_{w_k}$$
for $w=w_1w_2\cdots w_k\in W$.
It is known that a surjective coding map $$\tilde\pi:\Sigma\to K$$
is defined by $\tilde\pi(\omega_1\omega_2\cdots)=\lim_{k\to\infty}
g_{\omega_1\omega_2\cdots\omega_k}(z)$, which is independent of 
$z\in \tilde S$. 
\end{defi}

\begin{pro}\label{pro:ss}
$\phi\circ\tilde\pi=\pi$.
\end{pro}

\begin{proof}
We show that $\tilde l_w(1)=g_w(\tilde x)$ for any $w\in W$ 
by induction.
Assume that $\tilde l_w(1)=g_w(\tilde x)$ for $w\in W_k$.  
Then $\phi\circ g_w(\tilde x)=x_w$.
By (\ref{eq:0}), $f^k\circ \phi\circ g_w\circ \tilde l_i=l_i$.
Hence $F_{x_{w}}(l_i)=\phi\circ g_w\circ \tilde l_i$.
Therefore $l_{wi}=l_w\cdot(\phi\circ g_w\circ \tilde l_i)$ 
lifts to $\tilde l_w\cdot (g_w\circ \tilde l_i)$.
Hence $\tilde l_{wi}(1)=g_w\circ \tilde l_i(1)=g_{wi}(\tilde x)$.

From $\tilde \pi(w_1w_2\cdots)=\lim_{k\to\infty}g_{w_1w_2\cdots w_k}
(\tilde x)$ and $g_w(\tilde x)=\tilde l_w(1)$, 
it follows that $\phi\circ \tilde \pi(w_1w_2\cdots)
=\lim_{k\to\infty}x_{w_1w_2\dots w_k}=\pi(w_1w_2\cdots)$.
\end{proof}

\begin{rem}
Let $$[r]_N=(\tilde l_i(1))\in(\phi^{-1}f^{-1}(\bar x))^d.$$
 Since the iterated function system $\mathcal I$ is determined by 
$(\tilde l_i(1))_i$, 
so is the coding map $\pi_r$. 
\end{rem}

\section{Examples}\label{sec:3}

\subsection{}
Let $f(z)=z^d$.
Then $C=P=\mathrm{AP}=\{0,\infty\}$, $S=\mathbb C-\{0\}$, and 
the Julia set $J$ is the unit circle $\{|z|=1\}$.
Since $P-\mathrm{AP}=\emptyset$, we have a unique ramification function 
$\rho(x)=1$, and so $N=N^\rho$ is trivial.
Thus $\tilde S$ is Euclidean (i.e. $\tilde S=\mathbb C$).
We take an universal covering $\phi=z\mapsto e^{-2\pi iz}:\mathbb C\to 
\mathbb C-\{0\}$. 
Fix basepoints $\bar x=1\in\mathbb C-\{0\}$ and $\tilde x=0\in
\mathbb C$.
Then $\phi^{-1}f^{-1}(\bar x)=\frac1d\mathbb Z$.
The contraction associated with $n/d\in\frac1d\mathbb Z$ is 
$z\mapsto z/d+n/d$.
The coding maps $\pi_r$ such that $[r]_N=(n_1/d,n_2/d,\dots,n_d/d)$ with 
$\{n_1,n_2,\dots,n_d\}=\{0,1,\dots,d-1\}\mod d$, 
for which the attractors $K_r$ are intervals of length one, 
are considered canonical.

Consider the case $d=2$.
If $[r]_N=(n_1/2,n_2/2)$ with $n_1\le n_2$, 
then $K_r$ is the closed interval $[n_1,n_2]$, and $\pi_r$ 
is almost $(n_2-n_1)$-to-one map provided $n_2-n_1>0$.
Since $K_r$ is a closed interval, $K_r$ tiles $\phi^{-1}(J)=\mathbb R$ 
whenever $n_2-n_1$ is positive. 

Consider the case $d=3$.
Let $r$ be a radial with $[r]_N=(n_1/3,n_2/3,n_3/3)=(km_1/3,km_2/3,km_3/3)
$, where $k$ is the greatest common divisor of $n_1,n_2$, and $n_3$. 
Kenyon \cite{Ken97} proved that $\mu(K_r)$, 
the 1-dimensional Lebesgue measure of 
$K_r$ is positive if and only if $m_1+m_2+m_3=0\mod 3$, and that then 
$\mu(K_r)=k$ and $K_r$ tiles $\phi^{-1}(J)=\mathbb R$.
Remark that $K_r$ is not necessarily an interval.
Related topics are discussed by Lagarias and Wang \cite{LaWa96}.

\subsection{}
Let $f(z)=z^2-2$.
Then $C=\{0,\infty\}$, $P=\{-2,2,\infty\}$, $\mathrm{AP}=\{\infty\}$, 
$S=\mathbb C$, and 
the Julia set $J$ is the interval $[-2,2]$.
The canonical ramification function is $\rho(z)=2$ if $z=-2,2$, 
$\rho(z)=1$ otherwise.
Note that $\tilde S$ is Euclidean, 
and take a universal branched covering $\phi=z\mapsto 2\cos\pi z:
\mathbb C\to\mathbb C$.
Then $\phi^{-1}(J)$ is the real axis and $f\circ\phi(z)=\phi(\pm 2z+2n),
n\in\mathbb Z$.
Fix basepoints $\bar x=0$ and $\tilde x=1/2$.
Then $f^{-1}(\bar x)=\{\pm\sqrt{2}\}$, 
$\phi^{-1}(\sqrt{2})=\{1/4+2n,7/4+2n\,|\,n\in\mathbb Z\}$ 
and $\phi^{-1}(-\sqrt{2})=\{3/4+2n,5/4+2n\,|\,n\in\mathbb Z\}$.
The contraction associated with $\pm1/4+n\in \phi^{-1}f^{-1}(\bar x)$ 
is $z\mapsto \pm z/2+n$.
If $[r]_N=(k_1,k_2)\in\left(\phi^{-1}f^{-1}(\bar x)\right)^2$, 
then $K_r$ is the closed interval $[m_1,m_2]$, where
\begin{multline*}
m_1=\min\{2n_1,2n_2\},m_2=\max\{2n_1,2n_2\}
\text{ if }k_1=1/4+n_1,k_2=1/4+n_2, 
\end{multline*}
\begin{multline*}
m_1=\min\{2n_1,n_2-n_1\},m_2=\max\{2n_1,n_2-n_1\}\\
\text{ if }\{k_1,k_2\}=\{1/4+n_1,-1/4+n_2\},
\end{multline*}  
\begin{multline*}
 m_1=\min\{\frac{4n_1-2n_2}3,\frac{4n_2-2n_1}3\}, 
m_2=\max\{\frac{4n_1-2n_2}3,\frac{4n_2-2n_1}3\}\\
\text{ if }k_1=-1/4+n_1,k_2=-1/4+n_2.
\end{multline*}
The coding maps $\pi_r$ such that $[r]_N=(k_1,k_2)$ with $\{k_1,k_2\}
=\{1/4+n_1,-1/4+n_2\}$ ($|3n_1-n_2|=1$), for which the attractors $K_r$ 
are intervals of length one, are considered canonical. 

\subsection{}
\begin{figure}[hbtp]
\centering
\includegraphics[height=8cm]{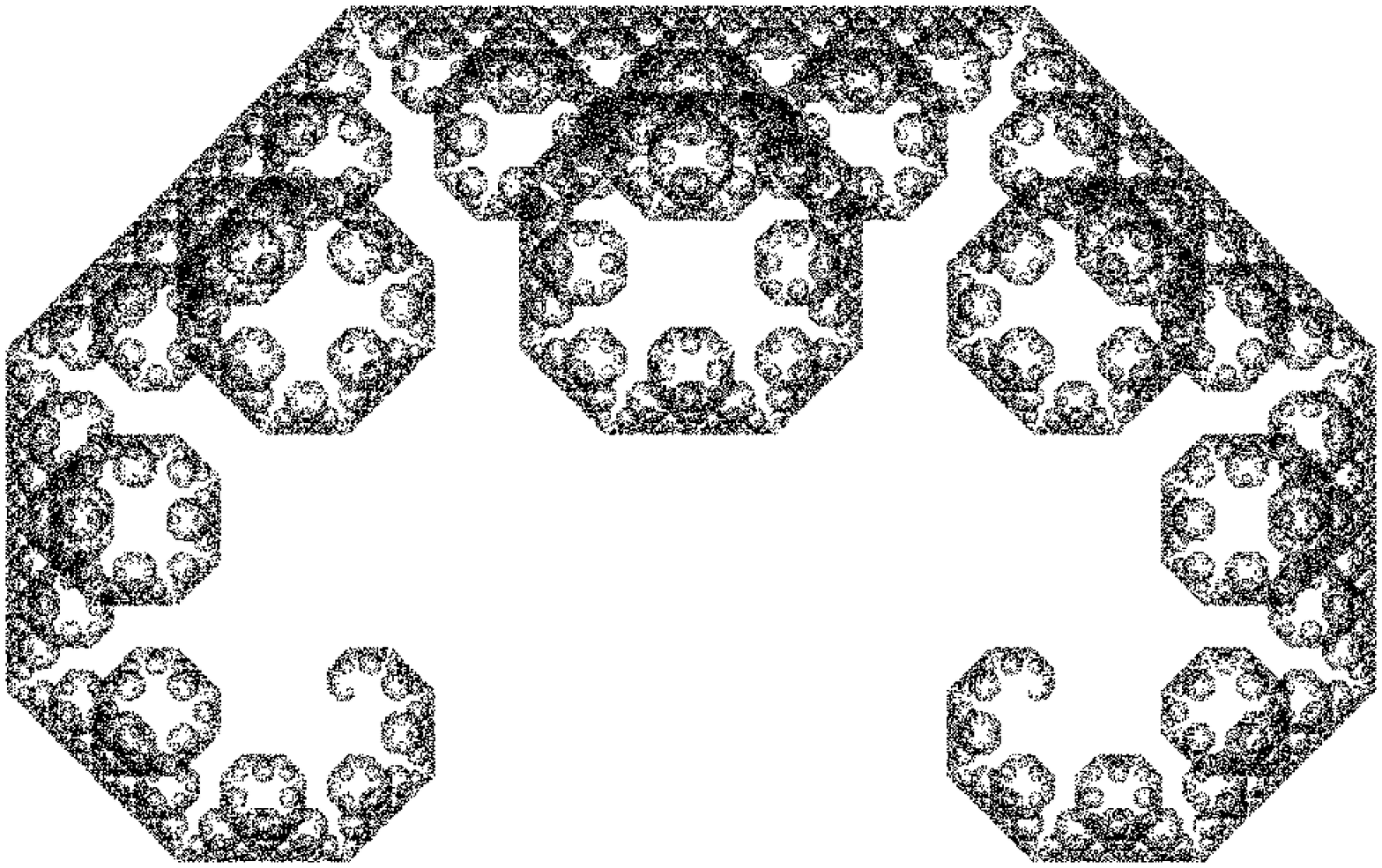}
\caption{The L\'evy Dragon.}
\label{fig:levy}
\end{figure}

Let $f(z)=-{(z-1)^2}/{4z}$.
(This map is conjugate to the map $z\mapsto (z-2)^2/z^2$, 
which is discussed in \cite{Be91}, \S4.3.)
Then $C=\{-1,1\},P=\{1,0,\infty\},\mathrm{AP}=\emptyset$ and 
$S=J=\hat{\mathbb C}$. 
The canonical ramification function is $\rho(0)=\rho(\infty)=4$, 
$\rho(1)=2$, and $\rho(z)=1$ otherwise.
Note that $\tilde S$ is Euclidean.
If $\phi:\mathbb C\to \hat{\mathbb C}$ is an elliptic function 
of order four with lattice $2\Gamma=\{2n+2mi\,|\,n,m\in\mathbb Z\}$ 
such that $\phi(iz)=\phi(z),
\phi(0)=\infty,\phi(1)=1,\phi(1+i)=0$, then $\phi(\alpha z+\beta)
=f(\phi(z))$ for $\alpha\in\{\pm1\pm i,\pm1\mp i\},
\beta\in2\Gamma$ and $\phi$ is a 
universal covering for the orbifold $(S,\rho)$.
Clearly, $\phi^{-1}(J)=\mathbb C$.
Note that $\phi(z)=a\wp(2z)^2$ satisfies the 
above condition, where $\wp$ is the 
Weierstrass elliptic function for the lattice $\Gamma=
\{n+mi\,|\,n,m\in\mathbb 
Z\}$ and $a$ is some constant.
It is easily seen that the map $f$ is a version of Latt\`es' example.
(For example, see \cite{Be91}, \S4.3 or \cite{Mi99}, \S7.
Using the addition formula, we have 
$\wp(\alpha z)^2=-(\wp(z)^2-g_2/4)^2/4\wp(z)^2$ for $\alpha\in
\{\pm1\pm i,\pm1\mp i\}$, where $g_2$ is a nonzero constant. 
The constant $a$ is equal to $4/g_2$.)
Fix basepoints $\bar x=-1$ and $\tilde x=1/2+i/2$.
Then $\phi^{-1}f^{-1}(\bar x)=\{s+n+mi\,|\,n,m\in\mathbb Z,s=1/2\text
{ or }i/2\}$.
The contraction associated with $x\in\phi^{-1}f^{-1}(\bar x)$ is 
$$\begin{array}{ll}
z\mapsto (1-i)z/2+n+mi& \text{ if }x=1/2+n+mi\\ 
z\mapsto (1+i)z/2+n+mi& \text{ if }x=i/2+n+mi\\
z\mapsto (-1+i)z/2+n+mi& \text{ if }x=-1/2+n+mi\\
z\mapsto (-1-i)z/2+n+mi& \text{ if }x=-i/2+n+mi,\\
\end{array}$$
where $n+mi\in(1+i)\Gamma=\{n+mi\,|\,n+m\text{ is even}\}$.
The attractor $K_r$ is a compact set with integral Lebesgue measure. 
For example, in the case $[r]_N=(-i/2+1+i,-1/2+2)$, $K_r$ is the triangle 
with vertices $0,2,1+i$.
This case is considered as canonical.
In the case $[r]_N=(i/2,1/2+1+i)$, $K_r$ is the L\'evy Dragon (Figure 
\ref{fig:levy}).
Thus Tiling Theorem  gives another proof of the well-known 
fact that the L\'evy Dragon tiles $\mathbb R^2$ and has nonempty interior.

We can calculate the 2-dimensional Lebesgue measure $\mu(K_r)$ for 
$[r]_N=(\alpha_1/2+\beta_1,\alpha_2/2+\beta_2)$ 
($\alpha_j\in\{\pm1,\pm i\},
\beta_j\in(1+i)\Gamma$) as follows. 
$$
\mu(K_r)=\left\{
\begin{array}{cl}
2|\beta_2-\beta_1|^2&\text{ if }\alpha_1=\alpha_2\\
\left|\frac{\beta_2}{2-(1-i)\alpha_2}
-\frac{\beta_1}{2-(1-i)\alpha_1}\right|^2
&\text{ if }(\alpha_1,\alpha_2)=(1,i),(i,1)\\
10\left|\frac{\beta_2}{2-(1-i)\alpha_2}
-\frac{\beta_1}{2-(1-i)\alpha_1}\right|^2
&\text{ if }(\alpha_1,\alpha_2)\in 
A_1\times A_2\text{ or }A_2\times A_1\\
25\left|\frac{\beta_2}{2-(1-i)\alpha_2}
-\frac{\beta_1}{2-(1-i)\alpha_1}\right|^2
&\text{ if }(\alpha_1,\alpha_2)=(-1,-i),(-i,-1)
\end{array}
\right.
$$
where $A_1=\{1,i\},A_2=\{-1,-i\}$.
The details are left to the reader.

\subsection{}\label{sec:3.1}
Let $f(z)=z^2-3$.
Then $C=\{0,\infty\}$, $P=\{-3,6,33,\cdots,\infty\}$, 
$\mathrm{AP}=\{\infty\}$, $S=\mathbb C$
,
 and the Julia set $J$ is a Cantor set in the real axis.
The canonical ramification function is $\rho(z)=2$ if $z=-3,6,33,\dots$,
$\rho(z)=1$ otherwise.
Thus $\tilde S$ is hyperbolic (i.e. $S=\{|z|<1\}$), and for any radial,
the corresponding contractions $g_1,g_2$ are not invertible.

\begin{figure}[hbtp]
\centering
\includegraphics[height=4cm]{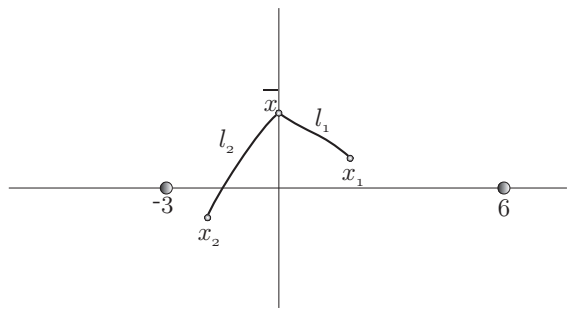}
\caption{The radial $r_1$.}
\label{fig:fig1}
\end{figure}

\begin{figure}[hbtp]
\centering
\includegraphics[height=4cm]{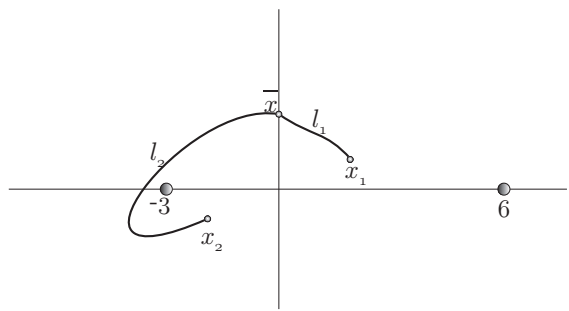}
\caption{The radial $r_2$.}
\label{fig:fig2}
\end{figure}

\begin{figure}[hbtp]
\centering
\includegraphics[height=4cm]{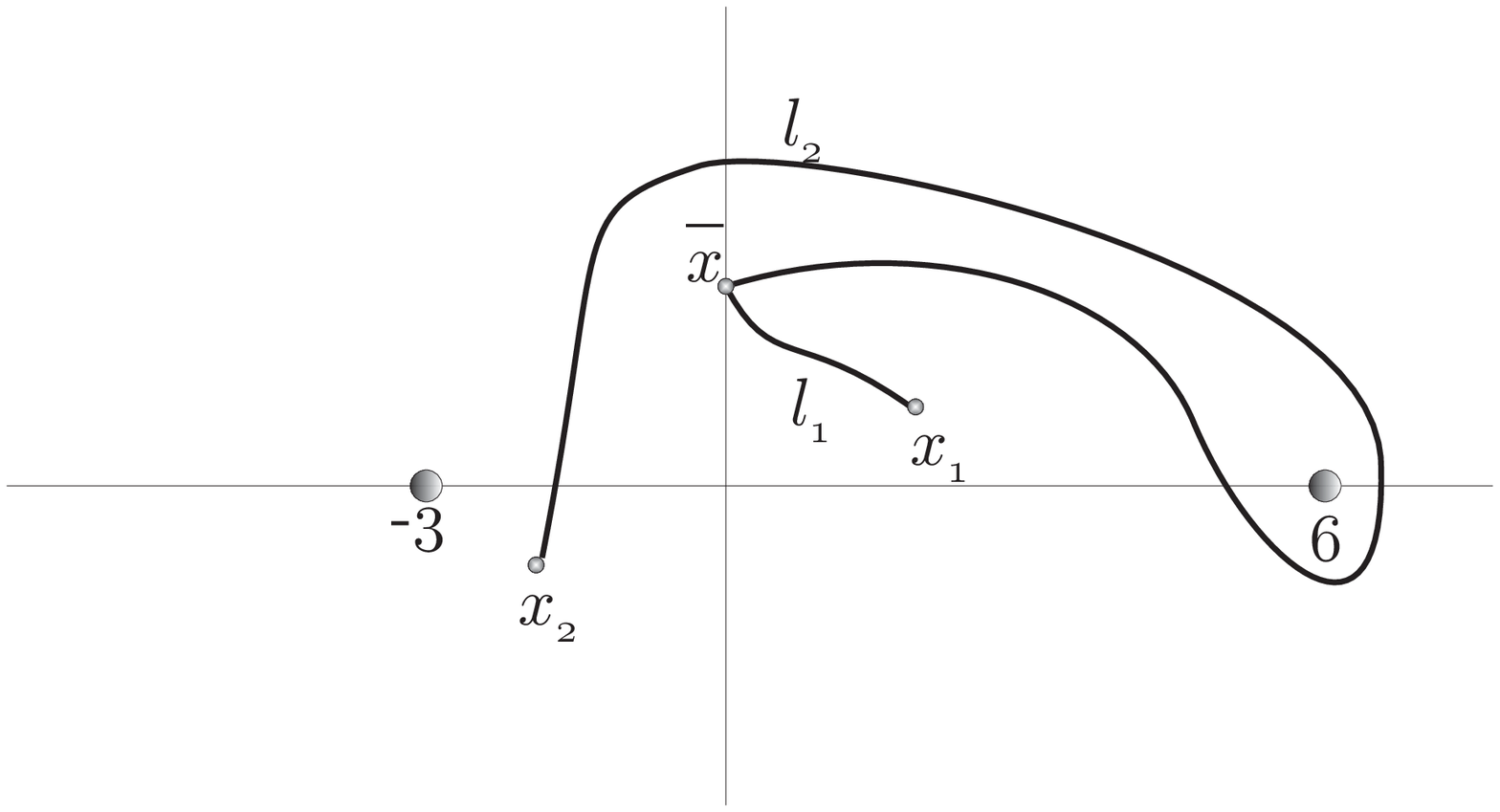}
\caption{The radial $r_3$.}
\label{fig:fig3}
\end{figure}

If we take radials $r_1$, $r_2$ and $r_3$ as in Figures \ref{fig:fig1}, 
\ref{fig:fig2} and \ref{fig:fig3}, then 
 $\pi_{r_1}$ is a homeomorphism, $\pi_{r_2}$ is exactly two-to-one, 
and $\pi_{r_3}$ is at most three-to-one respectively.
Moreover, $\#\pi_{r_3}^{-1}(x)=1$ for $\mu$-almost all $x$.
So, we say that $\pi_{r_3}$ is almost one-to-one.
See Section \ref{sec:5} for the proof.
The coding map $\pi_{r_1}$ is considered canonical.

\section{Further setting}\label{sec:4}

In this section, we define an invariant subgroup and 
the equilibrium state (the Brolin-Lyubich measure).

Let $f$ be a subhyperbolic rational map of degree $d$, and 
$\rho=\rho_f$ the canonical ramification function.
Fix a basepoint $\bar x$.

\begin{defi}
For a subgroup $N\subset G=\pi_1(S',\bar x)$, we have a covering 
$$\phi_N:(S'_N,x^N)\to(S',\bar x)$$
with $\phi_{N*}\pi_1(S'_N,x^N)=N$.
If $N^\rho\subset N$, then 
we can extend $\phi_N$ to a branched covering 
$$\phi_N:(S_N,x^N)\to (S,\bar x),$$
where $S_N-\phi_N^{-1}(P)=S_N'$.
For $l\in \bar L(f,\bar x)$, we denote by $l^N$ the lift of $l$ to $S'_N$ 
(or $S_N$) with $l^N(0)=x^N$.

If $N'$ is a subgroup with $N'\subset N\subset G$, 
then there exists a unique covering 
$$\psi_{N',N}:(S'_{N'},x^{N'})\to (S'_N,x^N)$$
such that $\phi_N\circ\psi_{N',N}=\phi_{N'}$.
If $N^\rho\subset N'$, then $\psi_{N',N}$ is extended to a branched 
covering $\psi_{N',N}:S_{N'}\to S_{N}$. 

Suppose $N$ is a normal subgroup.
The group of deck transformations for $\phi_N$ is identified with 
$G/N$.
For $\gamma\in\Gamma(f,\bar x)$, we denote the quotient class for 
$[\gamma]$ by $[\gamma]_N$.
\end{defi}

\begin{defi}
Let $l\in\Lambda(f,\bar x)$ with terminal point $x$.
Consider a homomorphism $f_*:\pi_1(S'-f^{-1}(P),x)
\to G$ induced from $f$ and a homomorphism 
$l_{\#}^{-1}\iota_*:\pi_1(S'-f^{-1}(P),x)
\to G$ induced from the inclusion and the path $l^{-1}$.
 A subgroup $N\subset G$ is said to be {\em invariant} with respect 
to $l$ if $N\subset f_*(l_{\#}^{-1}\iota_*)^{-1}(N)$, or 
equivalently if $F_x(\gamma)$ is a closed curve and 
 $[lF_x(\gamma)l^{-1}]\in N$ for every $\gamma$ with $[\gamma]\in N$.

If $N$ and $N'$ are invariant with respect to $l$, then so is 
the subgroup generated by $N$ and $N'$.
If $N$ and $N'$ are invariant with respect to $l$ and $l'$ respectively, 
then $N\cap N'$ is invariant both with respect to $l$ and $l'$.
We denote by $N_l$ the maximal invariant subgroup with respect to $l$.
A subgroup $N$ is invariant with respect to a radial $r=(l_i)$ if 
$N$ is invariant with respect to all $l_i,i=1,2,\dots,d$.
The maximal invariant subgroup $N_r$ with respect to $r$ is equal to 
$\bigcap_{i=1}^d N_{l_i}$. 
It is evident that $N^\rho$ is invariant with respect to every 
radial for any ramification function $\rho$.
\end{defi}

\begin{pro}\label{pro:12}
Let $l\in\Lambda(f,\bar x)$ and let $N\subset G$ be a subgroup.
If $N$ is invariant with respect to $l$, then there exists a 
covering $g:S_N'\to S'_N-\phi_N^{-1}f^{-1}(P)$ such that 
$f\circ\phi_N\circ g=\phi_N$ and  $g(x^N)=l_N(1)$.
If $N^\rho\subset N$, then $g$ is extended to a branched covering 
$g:S_N\to S_N$.
\end{pro}

\begin{proof}
Let $[\gamma]\in N$.
Then $[lF_x(\gamma)l^{-1}]\in N$, and so $lF_x(\gamma)l^{-1}$ 
can be lifted to a closed curve in $S'_N-\phi_N^{-1}f^{-1}(P)$. 
This means $(f\circ \phi_N)_*\pi_1(S'_N-\phi_N^{-1}f^{-1}(P),l_N(1))
\supset N=\phi_{N*}\pi_1(S'_N,x^N)$.
Therefore we have a covering $g:S'_N\to S'_N-\phi_N^{-1}
f^{-1}(P)$ with $f\circ\phi_N\circ g=\phi_N$ 
and  $g(x^N)=l_N(1)$. 
If $N^\rho\subset N$,  we can extend the map $g$ on $S_N$.
\end{proof}

\begin{rem}
The covering $g$ above is contractive in the pullback metric 
$\phi_N^*\nu$,  where $\nu$ is the expanding metric for $f$. 
\end{rem}

Let $r=(l_i)$ be a radial, and $N$ an invariant subgroup with 
respect to $r$.
There exist 
contractions $g_1,g_2,\dots,g_d:S_N'\to S_N'$ such that
$f\circ\phi_N\circ g_i=\phi_N$ and $g_i(x^N)=l_i^N(1)$.
Similarly to Proposition \ref{pro:ss}, 
\begin{equation}
 l_w^N(1)=g_w(x^N).
\label{eq:4-1}
\end{equation}
For simplicity, we assume $N^\rho\subset N$.

\begin{defi}
We have a nonempty compact set $K=K^N=K_r^N\subset S_N$ with
$K=\bigcup_{i=1}^d g_i(K)$, and a coding map $\pi^N=\pi_r^N:\Sigma\to K$ 
with $\pi^N(i\omega)=g_i\pi^N(\omega)$.
We say $K$ is a {\em Julia tile} if $\phi_N(K)=J$.
\end{defi}

From (\ref{eq:4-1}), we have 
\begin{equation*}
\pi=\phi_N\circ\pi^N.
\label{eq:4-2}
\end{equation*}

\begin{pro}\label{pro:13}
Let $r$ be a radial.
 The maximal invariant subgroup $N_r$ coincides with the normal 
subgroup 
$$N_r'=\{[\gamma]\in G\,|\, F_{x_w}(\gamma)\text{ is a closed curve for 
any }w\in W\}.$$
\end{pro}

\begin{proof}
 If $F_{x_w}(\gamma)$ is a closed curve for any $w\in W$, then 
$F_{x_w}(l_{x_i}F_{x_i}(\gamma)l_{x_i}^{-1})
=l_{x_{wi}}F_{x_{wi}}(\gamma)l_{x_{wi}}^{-1}$ is a closed curve 
for any $i=1,2,\dots,d$ and any $w\in W$.
Thus $N_r'$ is invariant with respect to $r$. 

If $[\gamma]\in N_r$, then $F_{x_i}(\gamma)$ is a closed curve 
and $[l_iF_{x_i}(\gamma)l_i^{-1}]\in N_r$ for $i=1,2,\dots,d$. 
Inductively, we see that $F_{x_w}(\gamma)$ is a closed curve and 
$[l_wF_{x_w}(\gamma)l_w^{-1}]\in N_r$ for $w\in W$. 
Thus $[\gamma]\in N_r'$.
\end{proof}

\begin{defi}
We consider the monodromy actions $\eta_k$ of $G$ on $f^{-k}(\bar x)$ for 
$k=1,2,\dots$, that is, $\eta_k([\gamma])(x)
=x\cdot [\gamma]=F_x(\gamma)(1)$.
We set 
$$\begin{array}{ccl}
\hat N&=&\bigcap_{k=1}^\infty\ker\eta_k\\
&=&\{[\gamma]\in G\,|\,\text{$F_x(\gamma)$ is a closed curve for any 
$k\ge1$ and any $x\in f^{-k}(\bar x)$}\}.
\end{array}$$
The quotient group 
$$\hat G=G/\hat N$$
is called the {\em reduced fundamental group} of $f$.
\end{defi}

From Proposition \ref{pro:13}, we immediately have

\begin{pro}
 If $l\in\Lambda(f,\bar x)$, then $\hat N$ is invariant with respect to 
$l$.
If $r=(l_i)$ is a proper radial, then $\hat N=N_r$.
\end{pro}

\begin{pro}\label{pro:1}
Let $l,l'\in\Lambda(f,\bar x)$, and let $N$ be an invariant subgroup 
with respect to $l$ and $l'$.
We have two coverings $g,g':S'_N\to S_N-\phi_N^{-1}f^{-1}(P)$ 
corresponding to $l,l'$ 
respectively by Proposition \ref{pro:12}. 
Suppose there exists a deck transformation $t:S'_N\to S'_N$ for the
 covering $\phi_N$ such that 
$t(z)=z'$ and $g(z)=g'(z')$ for some $z,z'\in  S'_N$.
Then  $g'\circ t=g$. 

In particular, if  $N$ is an invariant normal subgroup, then for the  
two coverings $g,g':S'_N\to S_N-\phi_N^{-1}f^{-1}(P)$ 
corresponding to $l,l'$, there exists a deck transformation $t$ 
with $g'\circ t=g$. 
\end{pro}

\begin{proof}
Write $y=g(z)=g'(z')$ and $x=\phi_N(z)=\phi_N(z')$. 
Consider the induced homomorphisms 
$$\begin{array}{rccc}
{g}_*:&\pi_1(S'_N,z)&\to &\pi_1(S'_N-\phi_N^{-1}f^{-1}(P),y),\\
{g'}_*:&\pi_1(S'_N,z')&\to&\pi_1(S'_N-\phi_N^{-1}f^{-1}(P),y),\\
 \phi_{N*}:&\pi_1(S'_N,z)&\to &\pi_1(S',x),\\
 \phi_{N*}':&\pi_1(S'_N,z')&\to &\pi_1(S',x),\\
(f\circ \phi_N)_*:&\pi_1(S'_N-\phi_N^{-1}f^{-1}(P),y)&\to &\pi_1(S',x).\\
\end{array}$$
The existence of $t$ implies $\phi_{N*}\pi_1(S'_N,z)=
\phi_{N*}'\pi_1(S'_N,z')$.
By (\ref{eq:0}) and the injectivity of $(f\circ\phi)_*$, we have 
 ${g}_*\pi_1(S'_N,z)={g}'_*\pi_1(S'_N,z')$.
Therefore there exists a unique homeomorphism $t':S_N\to S_N$
such that $g'\circ t=g$ and $t'(z)=z'$.
By the uniqueness of $t$, we have $t=t'$. 
\end{proof}

\begin{defi}
A Borel measure $\mu$ defined as follows is called the 
{\em Brolin-Lyubich measure} for $f$.
Let 
$$\mu_k=d^{-k}\sum_{y\in f^{-k}(\bar x)}\delta_y,$$
 where $\delta_y$ is the point mass at $y$. 
We define $\mu$ as the weak$*$ limit of $\mu_k$, that is, 
$\int hd\mu=\lim_{k\to\infty}\int hd\mu_k$ for any continuous function 
$h$.
It is known that the limit $\mu$ exists for any rational map of the 
Riemann sphere and $\mu$ is the unique equilibrium state, 
hence $f$ is strongly mixing with respect 
to $\mu$  (see \cite{Lj83}).
In our case, since $f$ is expanding, the direct proof is not difficult.
(We can also construct $\mu$ via a  Markov partition. 
Cf. \cite{Bo75}, \cite{Ja68}, \cite{Ja69}, \cite{Gu70}.)

\begin{pro}\label{pro:19}
The measure $\mu$ is characterized as the unique Borel probability 
measure satisfying 
\begin{equation}
\int h(x)d\mu(x)=d^{-1}\int \sum_{y\in f^{-1}(x)}h(y)d\mu(x)
\label{eq:equil}
\end{equation}
for any continuous function $h$.
\end{pro}

\begin{proof}
The operator $\mu\mapsto L(\mu)$ defined by 
$$\int h(x)dL(\mu)(x)=d^{-1}\int \sum_{y\in f^{-1}(x)}h(y)d\mu(x)$$ 
is a contraction of the space of Borel probability measures having  
compact supports in $S$ with metric $D(\mu,\mu')=\sup_{h}
|\int hd\mu-\int hd\mu'|$, where the supremum runs over all functions
 $h$ with Lipschitz constant (which is taken with respect to the
 expanding metric on $S$) 
less than or equal to one.
Indeed, if the Lipschitz constant of $h(x)$ is less than $\alpha$, then 
the Lipschitz constant of $d^{-1}\sum_{y\in f^{-1}(x)}h(y)$ is less than 
$c\alpha$, where $c$ is the contraction constant (see the end of 
Definition \ref{def:3}).
\end{proof}

\begin{rem}
Lyubich \cite{Lj83} proved that  
(\ref{eq:equil}) holds for any rational map $f$ 
and any Borel function $h(x)$.
\end{rem}

Let $N$ be a subgroup of $G$.
A Borel measure $\mu_N$ on $S'_N$ (or $S_N$) is defined as the lift of 
$\mu$, that is, for a small Borel set $E\subset S'_N$ for which the 
restriction $\phi|E$ is injective, we set $\mu_N(E)=\mu(\phi(E))$. 
We write 
$$J_N=\phi_N^{-1}(J).$$
It is easily seen that $\mu$ (or $\mu_N$) is supported on $J$ 
(or $J_N$).
\end{defi}

Hence we have

\begin{pro}\label{pro:3}\
\begin{enumerate}
 \item \label{pro:3-2}
$\mu$ is an invariant measure (i.e. $\mu(E)=\mu(f^{-1}(E))$ 
for any Borel set $E\subset S$). 
\item $f$ is strongly mixing with respect to $\mu$ (i.e.
$\lim_{k\to\infty}\mu(A\cap f^{-k}(B))=\mu(A)\mu(B)$ for Borel sets 
$A,B\subset S$).  
\item \label{pro:3-2.5}
For a Borel set $E\subset S$, $d\cdot\mu(E)=\mu(f(E))$ provided 
$f:E\to f(E)$ is injective.
\item \label{pro:3-3} 
Let $1\le i\le d$.
For a Borel set $E\subset S_N'$, 
$d\cdot\mu_N(g_i(E))=\mu_N(E)$ provided $g_i:E\to g_i(E)$ is 
injective. In general, if $g_i:E\to g_i(E)$ is at most $n$-to-one,
then $d\cdot\mu_N(g_i(E))\leq \mu_N(E)\le nd\cdot\mu_N(g_i(E))$.
\end{enumerate} 
\end{pro}

\begin{proof}
We need to show only (\ref{pro:3-2.5}).
(\ref{pro:3-3}) is an immediate consequence of (\ref{pro:3-2.5}).

Suppose $E$ is so small that there exists an open set $U\supset 
\overline E$ with $f:U\to f(U)$ injective.
Then for any continuous function $h$ with support included in $U$, 
we have $\int hd\mu=d^{-1}\cdot\int h\circ(f|U)^{-1}d\mu$  
by Proposition \ref{pro:19}.
Thus $\mu(E)=d^{-1}\mu(f(E))$.
In general, divide $E$ into small parts.
\end{proof}

\section{Julia tiling}\label{sec:tile}

In this section, we prove Tiling Theorem and Multiplicity Theorem.
Moreover, we give several necessary and sufficient conditions for 
suitableness of a coding map.

Let $f$ be a subhyperbolic rational map of degree $d$.
Let $r$ be a radial with basepoint $\bar x$, 
and $N$ an invariant normal subgroup with 
respect to $r$ such that $N^{\rho_f}\subset N$.
Write $K=K^N_r$.

\begin{thm}\label{thm:2}
 If $\pi(\Sigma)=J$, then there exists a subset $T$ of the group 
of deck transformations $G/N$ 
such that $\mu_N(t(K)\cap t'(K))=0$ for $t\ne t'$ and 
$\bigcup_{t\in T}t(K)=J_N$.
\end{thm}

\begin{proof}
By (\ref{pro:3-3}) of Proposition \ref{pro:3} 
and $K=\bigcup_{i=1}^dg_i(K)$, we have 
$$\mu_N(K)\le\sum_{i=1}^d\mu_N(g_i(K))\leq\mu_N(K).$$
Thus 
\begin{equation}
  \mu_N(K)=\sum_{i=1}^d\mu_N(g_i(K))=d\cdot\mu_N(g_i(K)),\ 
i=1,2,\dots,d,\label{eq:2}
\end{equation}
\begin{equation}
 \sum_{i\ne j}\mu_N(g_i(K)\cap g_j(K))=0,
\text{ and }
\sum_{i=1}^d\mu_N(\{x\in K\,|\,\#(g_i^{-1}(x)\cap K)\ge2\})=0.
\label{eq:3}\end{equation}

There exists a subset $U\subset K$ open 
in the relative topology of $J_N$. 
Indeed,  $\phi_N(K)=\phi_N(\pi^N(\Sigma))=\pi(\Sigma)=J$ implies 
$J_N=\bigcup_{t\in G/N}t(K)$.
Since $G/N$ is countable, $K$ has a nonempty interior in the relative 
topology by the Baire category theorem.
Take $\omega=\omega_1\omega_2\dots\in\Sigma$ with $\pi^N(\omega)\in U$.
If $k_0$ is sufficiently large, 
we have $g_w(K)\subset U$ for  $w=\omega_1\omega_2\cdots
\omega_{k_0}\in W_{k_0}$.  
Note that $\pi^N(w^\infty)\in U$. 
Therefore for any $x\in J_N$, there exists $k\ge0$ such that 
$g_w^k(x)\in U$.
Consequently, $\bigcup_{k=0}^\infty g_w^{-k}(K)=J_N$.

For $\mu_N$-almost all $x\in J_N$, there uniquely exists 
$y_x\in K$ such that $g_u(y_x)=g_w^k(x)$ for some $k\ge0$ and 
$u\in W_{kk_0}$.
Indeed, let 
$$E_k=\bigcup_{u\ne w\in W_k}(g_u(K)\cap g_w(K))\cup\bigcup_{w\in W_k}
\{x\in K\,|\, \#(g_w^{-1}(x)\cap K)\ge 2\}$$
and $E'=\bigcup_{w\in W}g_w^{-1}(\bigcup_{k>0}E_k)$.
Each $g_w:K\to K$ is at most finite-to-one.
Thus $\mu_N(E')=0$ by (\ref{pro:3-3}) of Proposition \ref{pro:3} 
and (\ref{eq:3}).
Let $x\in J_N-E'$ and let $k\ge0$ be the minimal integer such that
$g_w^k(x)\in K$. 
Since $g_w^k(x)\notin E_{kk_0}$, there  exists a unique $u\in W_{kk_0}$ 
such that 
$g_w^k(x)\in g_u(K)$ and $\#g_u^{-1}(g_w^k(x))\cap K=1$.
If $k'>k$ and $g_{u'}(y_x)=g_w^{k'}(x)$ for $u'\in W_{k'k_0}$, 
then necessarily $u'=w^{k'-k}u$, and so $g_u^{-1}(g_w^k(x))\cap K
=g_{u'}^{-1}(g_w^{k'}(x))\cap K$.

By Proposition \ref{pro:1}, there exists a unique deck transformation 
$t_x\in G/N$ such that $g_w^k\circ t_x=g_u$ and $t_x(y_x)=x$.
Let $T=\{t_x\,|\,x\in J_N-E'\}$.
Then $\bigcup_{t\in T}t(K)$ includes $J_N-E'$.
Since $\bigcup_{t\in T}t(K)$ is closed, it is equal to $J_N$.
By the uniqueness of $y_x$, we can see that $(t_x(K)\cap t_{x'}(K))-E'
=\emptyset$ whenever $t_x\ne t_{x'}$.
\end{proof}

\begin{co}
 If $\pi(\Sigma)=J$, then $K$ has nonempty interior in the relative 
topology of $J_N$.
\end{co}
\begin{proof}
 We have already shown this statement in the beginning of the second 
paragraph of the proof of Theorem \ref{thm:2}.
\end{proof}

\begin{co}\label{co:18}
 If $\pi(\Sigma)=J$, then $d\cdot\mu_N(g_i(A))=\mu_N(A),
1\le i\le d$ for any Borel set $A\subset K$.
\end{co}

\begin{proof}
Let $A\subset K$ be a Borel set.  
Let $E=E_1$ be the set defined in the proof of Theorem \ref{thm:2}.
Then $g_i:A-g_i^{-1}(E)\to g_i(A)-E$ is one-to-one.
Therefore $\mu_N(g_i(A))=\mu_N(g_i(A)-E)
=d^{-1}\mu_N(A-g_i^{-1}(E))=d^{-1}\mu_N(A)$ 
by (\ref{pro:3-3}) of Proposition \ref{pro:3}.
\end{proof}

\begin{defi}
 We denote by $m$ the {\em identically distributed Bernoulli measure} 
on $\Sigma$.
Namely, $m(\Sigma(w))=d^{-k}$ for every $w\in W_k$.
\end{defi}

\begin{pro}\label{pro:20}
If $\pi(\Sigma)=J$, then the conditional measure $\mu_N|K$ is 
the invariant measure for the iterated function system 
$\mathcal I$, that is, 
\begin{equation}
\mu_N(A\cap K)=d^{-1}\sum_{i=1}^d\mu_N(g_i^{-1}(A)\cap K)
\label{eq:4.3}
\end{equation}
for any Borel set $A$.
In particular, we have $\mu_N|K=\pi_{N*}m$ 
(i.e. $\mu_N(A)/\mu_N(K)=m(\pi_N^{-1}(A))$ 
for any Borel set $A\subset K$).
\end{pro}

\begin{proof}
(\ref{eq:4.3}) is an immediate consequence of Corollary \ref{co:18}.
It is easy to see that $\pi_{N*}m$ is also the invariant measure 
for $\mathcal I$ with weight $(1/d,1/d,\dots,1/d)$.
The uniqueness of the invariant measure (\cite{Hu81}) implies
 $\mu_N|K=\pi_{N*}m$.
\end{proof}

\begin{thm}\label{thm:3}
The following are equivalent.
\begin{enumerate}
 \item \label{it:thm3-0}
$\pi(\Sigma)=J$.
\item \label{it:thm3-0.5}
For distinct words $w,w'\in W$, $[l_w{l_{w'}}^{-1}]$ is nontrivial, and 
$\pi(\Sigma)\not\subset P$.
\item \label{it:thm3-1}
For distinct words $w,w'\in W$, $[l_w{l_{w'}}^{-1}]\notin N_r$.
\item \label{it:thm3-3}
$\lim_{k\to\infty}(\#\{g_w(x^N)\,|\,w\in W_k\})^{1/k}=d$ 
for $N=\{1\}$, and $\pi(\Sigma)\not\subset P$.
\item \label{it:thm3-3.5}
$\lim_{k\to\infty}(\#\{g_w(x^N)\,|\,w\in W_k\})^{1/k}=d$ 
for $N=N_r$.
\item \label{it:thm3-4}
For any $x\in J$, $\pi^{-1}(x)$ is finite.
\item \label{it:thm3-5}
$m(\{\omega\in\Sigma\,|\,\pi^{-1}\pi(\omega)\text{ 
is at most countable}\})>0$.
\item \label{it:thm3-6}
$\mu(\pi(\Sigma))>0$.
\end{enumerate} 
\end{thm}

\begin{proof}
 The implications (\ref{it:thm3-0})$\Rightarrow$(\ref{it:thm3-6})
and (\ref{it:thm3-4})$\Rightarrow$(\ref{it:thm3-5}) 
are trivial.

(\ref{it:thm3-0.5})$\Rightarrow$(\ref{it:thm3-3}),
(\ref{it:thm3-1})$\Rightarrow$(\ref{it:thm3-3.5}):
(\ref{it:thm3-0.5}) or (\ref{it:thm3-1}) implies 
$\#\{g_w(x^N)\,|\,w\in W_k\}=d^k$. 
Recall (\ref{eq:4-1}).

(\ref{it:thm3-3})$\Rightarrow$(\ref{it:thm3-0.5}),
(\ref{it:thm3-3.5})$\Rightarrow$(\ref{it:thm3-1}):
If $[l_wl_{w'}^{-1}]\in N$ for some $w,w'\in W_n$, then $g_w=g_{w'}$.
Thus $\#\{g_u(x^N)\,|\,u\in W_{kn}\}\le (d^n-1)^{k}$.

(\ref{it:thm3-6})$\Rightarrow$(\ref{it:thm3-0}):
Since $f\pi(\Sigma)=\pi(\Sigma)$, we have $\mu(\pi(\Sigma))=1$ 
by the ergodicity.
Thus (\ref{it:thm3-0}) follows from the compactness of $\pi(\Sigma)$.

(\ref{it:thm3-1})$\Rightarrow$(\ref{it:thm3-6}):
Since $|l_w|$ is bounded for $w\in W$, it follows from 
(\ref{it:thm3-1}) that there exists $n>0$ such that 
\begin{equation}
\#\{w'\,|\,
x_{w'}=x_w\}\le n
\label{eq:n}
\end{equation}
 for any $w\in W$.

Let $V$ be an open set including $\pi(\Sigma)$.
Take an open set $V'$ such that $\pi(\Sigma)\subset V'$ and 
$\overline {V'}\subset V$.
Then there exists $k_0>0$ such that if $k>k_0$, then $x_w\in 
V'$ for every $w\in W_k $.
Hence $$d^k/n\le\#\{x_w\,|\,w\in W_k\}\cap V'\le\#f^{-1}(\bar x)\cap V'.$$
Therefore $1/n\leq\mu(V)$, and so $1/n\leq\mu(\pi(\Sigma))$.

(\ref{it:thm3-1})$\Rightarrow$(\ref{it:thm3-0.5}):
Since (\ref{it:thm3-1})$\Rightarrow$(\ref{it:thm3-0}), 
(\ref{it:thm3-1}) implies $\pi(\Sigma)\not\subset P$.

(\ref{it:thm3-5})$\Rightarrow$(\ref{it:thm3-1}):
Suppose there exist distinct $w,w'\in W$ such that $[l_w{l_{w'}}^{-1}]
\in N_r$.
Set 
$$\Omega=\{\omega\in\Sigma\,|\,\sigma^k\omega\in \Sigma(w)\text{ 
for inifinitely many $k$}\}.$$
Then $m(\Omega)=1$.
For any $\omega=u_1wu_2wu_3\cdots\in\Omega$, 
$\pi^{-1}\pi(\omega)$ includes an uncountable set 
$\{u_1w_1u_2w_2u_3\cdots\,|\,w_i=w\text{ or }w'\}$.
Indeed, note that if $[l_v{l_{v'}}^{-1}],[l_u{l_{u'}}^{-1}]\in  N_r$, 
then $[l_{vu}{l_{v'u'}}^{-1}]\in N_r$.
Thus $x_{u_1wu_2w\cdots u_kw}=x_{u_1w_1u_2w_2\cdots u_kw_k}$ if 
$w_i=w\text{ or } w'$. 
Hence $\pi(\{u_1w_1u_2w_2u_3\cdots\,|\,w_i=w\text{ or }w'\})=\{
\pi(\omega)\}$.

(\ref{it:thm3-0.5})$\Rightarrow$(\ref{it:thm3-1}):
Suppose $\pi(\Sigma)\not\subset P$ and 
there exist distinct $w,w'\in W$ such that $[l_w{l_{w'}}^{-1}]
\in  N_r$.
If we take a large $k_0$, then $|F_{x_u}(l_{w}{l_{w'}}^{-1})|$ is so 
small for $u\in W_k,k\ge k_0$ 
that $F_{x_u}(l_{w}{l_{w'}}^{-1})$ is either homotopically trivial or 
 winding around some point of $P$ several times in $S'$.
Since $\pi(\Sigma)\not\subset P$, there exists a word $u$ such that 
$F_{x_u}(l_{w}{l_{w'}}^{-1})$ is homotopically trivial.
Then $[l_{uw}{l_{uw'}}^{-1}]$ is trivial.

(\ref{it:thm3-6})$\Rightarrow$(\ref{it:thm3-1}):
Suppose there exists distinct $w,w'\in W$ such that $[l_w{l_{w'}}^{-1}]
\in N_r$.
Then $g_w=g_{w'}$ by (\ref{eq:4-1}). 
From (\ref{eq:3}), we have  $\mu_N(g_w(K))=0$, and hence $\mu_N(K)=0$ 
by (\ref{eq:2}).
Therefore $\mu(\pi(\Sigma))=0$.

(\ref{it:thm3-1})$\Rightarrow$(\ref{it:thm3-4}):
Let $p\in J$.
Since $|l_\omega|$ is bounded for $\omega\in\Sigma$, 
$L(p)=\{l_\omega\,|\,\omega\in\pi^{-1}(p)\}$ is divided into 
a bounded number of  homotopy classes mod $N_r$, that is, 
$$L(p)=L(p)_1\sqcup L(p)_2\sqcup\cdots \sqcup L(p)_{m(p)}\quad
\left(l_\omega,l_{\omega'}\in L(p)_i\iff [l_\omega{l_{\omega'}}^{-1}]
\in N_r\right)$$
with $m(p)\le M$ for some $M>0$ independent of $p$.

Set
$$\begin{array}{ccl}
 A(p,k)&=&
  \{w\in W_k\,|\,\Sigma(w)\cap\pi^{-1}(p)\ne\emptyset\}\\
&=&\{w\in W_k\,|\, \text{ there exists $\omega\in \pi^{-1}f^k(p)$ 
such that $\pi(w\omega)=p$}\},\\
B(p,k,\omega)&=&\{x_w\,|\,w\in W_k,\pi(w\omega)=p\}.
\end{array}
$$
Then $\#B(p,k,\omega)$ is equal to or less than $\deg_pf^k$, 
the local degree of $f^k$ at $p$.
Indeed, if $\pi(\omega)=f^k(p)$, then we have exactly $\deg_pf^k$ lifts 
of $l_\omega$ by $f^{k}$ with terminal point $p$.
Thus there exists $b>0$ independent of $p$, $k$ and $\omega$ 
such that $\#B(p,k,\omega)<b$.
It is easily seen that $B(p,k,\omega)=B(p,k,\omega')$ whenever 
$[l_\omega{l_{\omega'}}^{-1}]\in N_r$.
Therefore choosing $\omega_i$ so that $l_{\omega_i}\in L(f^k(p))_i$, 
we have 
$$\begin{array}{ccl}
\#A(p,k)&\le &n\cdot\#\{x_w\,|\,w\in A(p,k)\}\\
&=&n\cdot\#\bigcup_{\omega\in\pi^{-1}f^k(p)}B(p,k,\omega)\\
&=&n\cdot\#\bigcup_{i=1}^{m(f^k(p))}B(p,k,\omega_i)\\
&\le&nMb.
\end{array}$$
for any $k>0$, where $n$ is defined in (\ref{eq:n}).
Hence 
\begin{equation}\#\pi^{-1}(p)\le nMb,
\label{eq:l}
\end{equation}
since $\pi^{-1}(p)= \bigcap_{k=1}^\infty
\bigcup_{w\in A(p,k)}\Sigma(w)$.
\end{proof}

\begin{thm}\label{thm:21}
 There exists an integer $n\ge 0$ such that for any invariant subgroup 
$N\subset G$ with respect to $r$,
\begin{enumerate}
\item $\#\pi^{-1}(p)=\#\phi_N^{-1}(p)\cap K=n$ for $\mu$-almost 
all $p\in J$, 
\item $\#\pi^{-1}(p)\ge\#\phi_N^{-1}(p)\cap K\ge n$ 
for any $p\in J-P$, and
\item $\mu_N(\phi_N^{-1}(A)\cap K)=n\mu(A)$ for any Borel set $A\subset 
\hat{\mathbb C}$.
In particular, $\mu_N(K)=n$.
\end{enumerate}
\end{thm}

\begin{proof}
First remark that $\pi=\phi_N\circ\pi^N$ and $\phi_N=\psi_{N,N_r}
\circ\phi_{N_r}$ imply $\#\pi^{-1}(p)\ge\#
\phi_N^{-1}(p)\cap K^N\ge\#\phi_{N_r}^{-1}(p)\cap K^{N_r}$ 
for every $p\in J$.
Let $E_k\subset S_{N_r}$ be the set as in the proof of Theorem 
\ref{thm:2}.
Then $D=\bigcup_{k=1}^\infty E_k$ is a null set.
If $x\in K^{N_r}-D$, then $\#({\pi^{N_r}})^{-1}(x)=1$.
Therefore $\#\pi^{-1}(p)=
\#\phi_{N_r}^{-1}(p)\cap K^{N_r}$ for every $p\in J-\phi_{N_r}(D)$.

If $\pi(\Sigma)\ne J$, then $\mu(\pi(\Sigma))=0$ by Theorem 
\ref{thm:3}, so the conditions are satisfied for $n=0$.

Suppose $\pi(\Sigma)=J$.
We show that  $h(x)=\#\phi_{N_r}^{-1}(x)\cap K$ is a Borel function.
To this end, let $A_{k,\epsilon}$ be the set of $x\in J$ 
such that for some $z_i\in K,i=1,2,\dots,k$, 
$\phi_{N_r}(z_i)=x,1\le i\le k$ and
the distance between $z_i$ and $z_j$ is  equal to or bigger than 
$\epsilon$ for $0\le i\ne j\le k$.
Then  $A_{k,\epsilon}$ is closed, and so  
 $\{x\,|\,h(x)\ge k\}=\bigcup_{n=1}^\infty A_{k-1,1/n}$ 
is a Borel set for $k\ge 2$.

For a Borel set $E\subset J$,
using $h(x)=d^{-1}\sum_{y\in f^{-1}(x)}h(y)$ and $\int_Eh(x)d\mu(x)=
\int_{f^{-1}(E)}h\circ f(y)d\mu(y)$ (substitute $h$ of (\ref{eq:equil}) 
by $1_{f^{-1}(E)}\cdot h\circ f$), we have
$$\int_Eh(x)d\mu(x) = \int_{f^{-1}(E)}h(y)d\mu(y).$$
Hence for any $E$ and any $k>0$, we have 
$\int_Eh(x)d\mu(x)=\int_{f^{-k}(E)}d\mu(x)$, 
which converges to $\mu(E)\int_Jh(x)d\mu(x)$ as $k\to\infty$ 
by the strong mixing condition.
Therefore $h(x)$ is constant almost everywhere.
Thus (1) is verified.
By the definition of $\mu_N$, we see that $\mu_N(\phi_N^{-1}(A)
\cap K)=n\mu(A)$ for any Borel set $A$.

If $z_i,z_i'\in K,i=1,2,\dots$ satisfy $\lim_{i\to\infty}z_i=
\lim_{i\to\infty}z_i'=z$, $z_i\ne z_i'$, and $\phi_N(z_i)=\phi_N(z_i')$, 
then $\phi_N(z)\in P$.
Since $\{x\,|\,\#\phi_N^{-1}(x)\cap K=n\}$ is dense in $J$, we have 
$\#\phi_N^{-1}(x)\cap K\ge n$ for every $x\in J-P$.
In general, $\#\phi_N^{-1}(x)\cap K\ge \max\{n-\max_{y\in \phi_N^{-1}(x)}
\deg_{y}\phi_N+1,1\}$ for every $x\in J$.
\end{proof}

\begin{defi}
 We call the number $n$ above the {\em multiplicity} of $\pi_r$ 
and denote by $n_r$ or $\mathrm{mul}(\pi_r)$.
\end{defi}

\begin{co}
 If $\pi(\Sigma)=J$, then $\mu=\pi_*(m)$.
\end{co}

\begin{proof}
By (3) of Theorem \ref{thm:21} and Proposition \ref{pro:20},
$m(\pi^{-1}(A))=\tilde\mu(\phi^{-1}(A)\cap K)/{n_r}=\mu(A)$ for 
any Borel set $A\subset\hat{\mathbb C}$. 
\end{proof}

\begin{pro}
If $\pi(\Sigma)=J$, then  the multiplicity of $\pi_r$ is equal to or 
bigger than $$\max_{x\in\{x_1,x_2,\dots,x_d\}}\#\{i\,|\,x=x_i\}.$$
In particular, if $r$ is not proper, then $n_r\ne 1$.
\end{pro}

\begin{proof}
Suppose $x_1=x_2=\cdots =x_s$.
Then $\pi(i\omega)=\pi(j\omega)$ for any $\omega\in\Sigma$ whenever 
$i,j\in\{1,2,\dots,s\}$.
Thus $\#\pi^{-1}\pi(i\omega)\ge s$ for $1\le i\le s$ and $\omega\in
\Sigma$.
Since $\mu(\pi(\Sigma(i)))\ge m(\Sigma(i))=d^{-1}>0$, 
the multiplicity is equal to or bigger than $s$.
\end{proof}

\section{Structure of $\mathrm{Cod}(f)$}\label{sec:st}

In this section, we prove Structure Theorem, and show that the monoid 
$\mathcal A(f)$ naturally acts on $\mathrm{Cod}(f)$.
The structures of $\mathrm{Cod}(f)$ for a couple of examples are 
closely investigated.

\begin{defi}
Let $N\subset G$ be a subgroup.
We say that $l$ and $l'\in \Lambda(f,\bar x)$ are 
{\em homotopic modulo $N$ with basepoint held fixed} 
if $l(1)=l'(1)$ and $[l_i{l_i'}^{-1}]\in N$.
We denote by $\Lambda_N(f,\bar x)$ the set of homotopy classes 
 modulo $N$ with basepoint $\bar x$ held fixed.
Then $\Lambda_N(f,\bar x)$ is identified with 
$\phi_N^{-1}f^{-1}(\bar x)$, and with $\{g:S'_N\to S'_N\,|\,
f\circ\phi_N\circ g=\phi_N\}$ as well.
The natural projection from $\Lambda(f,\bar x)$ to $\Lambda_N(f,\bar x)$ 
is denoted by 
$$l\mapsto [l]_N.$$

We say that two radials $r=(l_i)$ and $r'=(l_i')\in\mathrm{Rad}
(f,\bar x)=\Lambda(f,\bar x)^d$ are 
{\em  homotopic modulo $N$ with basepoint held fixed} 
if $[l_i]_N=[l_i']_N$ for $i=1,2,\dots,d$.
We write $$\mathrm{Rad}_N(f,\bar x)=\Lambda_N(f,\bar x)^d\text{ and }
[r]_N=([l_i]_N)_i.$$
\end{defi}

\begin{defi}\label{def:6}
 We say that two radials $r\in\mathrm{Rad}(f,\bar x)$ and 
$r'\in\mathrm{Rad}
(f,\bar x')$ are {\em freely homotopic} if there exists a 
homotopy 
$H:Q\times[0,1]\to S'$ between $r$ and $r'$ such that $H(\cdot,t):
Q\to S'$ is a radial for every $0\le t\le1$. 
\end{defi}

Let 
$$Y=\{\gamma:[0,1]\to S'\,|\,\gamma(0)=\bar x,\gamma(1)=\bar x'\}
\subset L.$$
We define an operation of $\gamma\in Y$ from $\Lambda(f,\bar x)$ to 
$\Lambda(f,\bar x')$ by 
$$\gamma\cdot l=\gamma^{-1} l F_{x_i}(\gamma).$$ 
For a subgroup $N\subset G$, the operation $\gamma:\Lambda(f,\bar x)
\to\Lambda(f,\bar x')$ naturally descends to 
$\gamma:\Lambda_N(f,\bar x)\to\Lambda_{\gamma_{\#}^{-1}(N)}(f,\bar x')$. 
Suppose $N'$ is invariant with respect to $l$ and $\gamma,\gamma'\in Y$ 
satisfy $[\gamma{\gamma'}^{-1}]\in N'$.
Then for $N\supset N'$, $\gamma\cdot[l]_{N}=\gamma'\cdot[l]_{N}$.
In particular, in the case $\bar x=\bar x'$,  the operation 
\begin{equation}
[\gamma]_{\hat N}:\Lambda_{N}(f,\bar x)\to\Lambda_{\gamma^{-1}_{\#}
(N)}(f,\bar x)\label{eq:6.1}
\end{equation} 
is well-defined for every $[\gamma]_{\hat N}\in \hat G=G/\hat N$.
Identifying $\hat G$ with the group of deck transformations of 
$\phi_{\hat N}:S_{\hat N}\to S$, we write the 
action $t:\Lambda_N(f,\bar x)\to\Lambda_{tN}(f,\bar x)$  by 
$$t\cdot g_i=tg_it^{-1}$$
for $t\in \hat G$, 
where we use the identification 
$\Lambda_N(f,\bar x)=\{g:S_N\to S_N\,|\,f\circ\phi_N\circ g=\phi_N\}$.
Clearly, the operation of $\gamma\in Y$ from $\mathrm{Rad}(f,\bar x)$ to 
$\mathrm{Rad}(f,\bar x')$ is defined diagonally.
From (\ref{eq:6.1})
\begin{equation}
[\gamma]_{\hat N}:\{[r]_{N_r}\,|\,r\in\mathrm{Rad}(f,\bar x)\}\to
\{[r]_{N_r}\,|\,r\in\mathrm{Rad}(f,\bar x)\}
\label{eq:6.2}
\end{equation} 
is well-defined for every $[\gamma]_{\hat N}\in \hat G$.

Two radials $r\in\mathrm{Rad}(f,\bar x)$ and $r'\in\mathrm{Rad}
(f,\bar x')$ are freely homotopic  if and only if 
there exists $\gamma\in Y$ such that $\gamma\cdot [r]_e=[r']_{e'}$, where 
$e$ and $e'$ are the trivial subgroups.
Thus we have a generalization of Definition \ref{def:6}.

\begin{defi}\label{def:19}
Let $N\subset G$ be a subgroup.
We say that $r$ and $r'$ are {\em freely homotopic modulo $N$} 
if there exists 
$\gamma\in Y$ such that $\gamma\cdot[r]_N=[r']_{\gamma_{\#}^{-1}(N)}$.
\end{defi}

\begin{thm}\label{thm:24}
Let $r=(l_i)\in\mathrm{Rad}(f,\bar x)$ and $r'=(l_i')\in
\mathrm{Rad}(f,\bar x')$.
Then $\pi_r=\pi_{r'}$ if and only if $r$ and $r'$ are freely homotopic 
modulo $N_r$. 
\end{thm}

\begin{proof}
For $r'$, the notation $x_w'$, $l'_w$, $l'_\omega$ and $F'_x(\cdot)$ 
are defined in a trivial way.

Suppose $\pi=\pi_r=\pi_{r'}$.
For $\omega\in\Sigma$, write $\gamma_\omega=l_\omega{l'_\omega}^{-1}$. 
Remark that 
\begin{equation}
F_{x_w}(\gamma_\omega)=F_{x_w}(l_\omega)
F'_{x'_w}(l'_\omega)^{-1}
\label{eq:6-2}
\end{equation}
 for $w\in W$ and $\omega\in\Sigma$ with 
$\pi(\omega)\notin P$, and 
\begin{equation}
l_wF_{x_w}(l_\omega)F'_{x'_w}(l'_\omega)^{-1}{l'_w}^{-1}
=\gamma_{w\omega}
\label{eq:6-3}
\end{equation}
 for $w\in W$ and $\omega\in\Sigma$.
Since the curve $F_{x_w}(l_\omega)F'_{x'_w}(l'_\omega)^{-1}$ 
joins $x_w$ and $x'_w$, 
the curve $F_{x_w}(l_\omega)F'_{x'_w}(l'_\omega)^{-1}F'_{x'_w}
(l'_{\omega'})F_{x_w}(l_{\omega'})^{-1}$ is closed for 
any $w\in W$ and $\omega,\omega'\in\Sigma$.
It follows from this that 
$[\gamma_\omega{\gamma_{\omega'}}^{-1}]\in N_r$ provided 
$\pi(\omega),\pi(\omega')\notin P$.

First we assume $\pi(\Sigma)\not\subset P$.
By (\ref{eq:6-2}) and (\ref{eq:6-3}),
 $\gamma_\omega^{-1}l_iF_{x_i}(\gamma_\omega)$ and 
$\gamma_\omega^{-1}\gamma_{i\omega}{l_i'}$ are homotopic in $S'$ 
with endpoints held fixed for $i=1,2,\dots,d$ 
whenever $\pi(\omega)\notin P$. 
Hence, since $[\gamma_\omega^{-1}\gamma_{i\omega}]\in 
\gamma^{-1}_{\omega\#}N_r$, we have
 $\gamma_\omega\cdot[r]_{N_r}=[r']_{\gamma_{\omega\#}^{-1}(N_r)}$.
If $\pi(\Sigma)\subset P$, then $\pi(\Sigma)$ consists of one fixed
 point (see Remark \ref{rem:5}), say $p$.
Note that $f$ is one-to-one near $p$.
Therefore if we modify $\gamma_\omega$ near $p$ into $\gamma'_\omega$ 
avoiding $P$, 
then $F_{x_w}(\gamma_\omega')$ coincides with $F_{x_w}(l_\omega)F'_{x'_w}
(l'_\omega)^{-1}$ 
except near $p$.
Hence we use the same discussion as above to conclude that $r$ and $r'$ 
are freely homotopic modulo $N_r$.

Suppose $r$ and $r'$ are freely homotopic modulo $N_r$.
Then there exists $\gamma\in Y$ such that 
$[\gamma l'_iF_{x_i}(\gamma)^{-1}
{l_i}^{-1}]\in {N_r}$ for $i=1,2,\dots,d$.
Let $x'\in\phi_{N_r}^{-1}(\bar x')$ be the terminal point of 
the lift $\gamma^{N_r}$ of $\gamma$ to $S_{N_r}$ whose initial point 
is $x^N$.
We have the contractions $g'_1,g'_2,\dots,g'_d:S_{N_r}\to S_{N_r}$  
corresponding to $r'$ with respect to the basepoints 
$\bar x'$ and $x'$.
Let $l_i^{'N_r}$ denote the lift of $l_i'$ to $S_{N_r}$ whose initial 
point is $x'$, and $l_i^{N_r}\gamma_0^{N_r}$ the lift of 
$l_iF_{x_i}(\gamma)$ to $S_{N_r}$ whose initial point is $x^{N_r}$.
Then the terminal points of $l_i^{N_r},l_i^{'N_r}$ are 
$g_i(x^{N_r}),g_i'(x')$ respectively.
By assumption, the terminal point of $\gamma^{N_r}l_i^{'N_r}$ coincides 
with that of $l_i^{N_r}\gamma_0^{N_r}$.
Thus $\gamma_0^{N_r}$ joins $g_i(x^{N_r})$ and $g_i'(x')$.
Since $\gamma_0^{N_r}$ is a lift of $\gamma$ by $f\circ\phi_{N_r}$, 
we have $g_i(\gamma^{N_r})=\gamma_0^{N_r}=g_i'(\gamma^{N_r})$.
Hence $g_i=g_i'$.
Consequently, $\pi_r^{N_r}=\pi_{r'}^{N_r}$, and so 
$\pi_r=\phi_{N_r}\circ\pi_r^{N_r}=\phi_{N_r}\circ\pi_{r'}^{N_r}=\pi_{r'}$.
\end{proof}

\begin{co}\label{co:34-1}
$$\mathrm{Cod}(f)\approx\{[r]_{N_r}\,|\,
r\in\mathrm{Rad}(f,\bar x)\}/\hat G,$$
where the action of $\hat G$ is defined in (\ref{eq:6.2}).
\end{co}

\begin{co}\label{co:35}
 If $\pi_r=\pi_{r'}$, then $N_r$ and $N_{r'}$ are conjugate 
(i.e. $N_{r'}=\gamma^{-1}_{\#}(N_r)$ for some $\gamma\in Y$).
\end{co}

\begin{exa}\ 
\begin{enumerate}
 \item $f(z)=z^d$.
Since $N^{\rho_f}$ is trivial, the fundamental group $G^{\rho_f}$ is equal
to $G=\mathbb Z$.
For any curve $l\in\Lambda(f,\bar x)$, the maximal invariant subgroup 
$N_l$ is trivial.
Thus $\hat G=G^{\rho_f}$ and 
$\mathrm{Cod}(f)\approx (\phi^{-1}f^{-1}(\bar x))^d/\hat G=
\mathbb Z^d/\mathbb Z$, where $n/d\in\phi^{-1}f^{-1}(\bar x)$ is 
identified with $n\in\mathbb Z$, and $\mathbb Z$ acts on $\mathbb Z^d$ by 
$$n\cdot(n_1,n_2,\dots,n_d)=(n_1+(d-1)n,n_2+(d-1)n,\dots,n_d+(d-1)n).$$
Thus $$\mathrm{Cod}(f)\approx\{[n_1,n_2,\dots,n_d]\,:\, n_i\in\mathbb Z,
0\le n_1\le d-2\},$$ where $[n_1,\dots]$ denotes the equivalence class 
including $(n_1,\dots)$.

\item
$f(z)=z^2-2$.
The fundamental group $G^{\rho_f}$ is equal to $\mathrm{Iso}(2\mathbb Z)=
\{x\mapsto ax+2n\,|\,a=\pm1,n\in\mathbb Z\}$.
For $l\in\Lambda(f,\bar x)$, 
let  $l_\infty=\phi(x)$ be the fixed point of $f$ such that 
$g(x)=x$ for the contraction $g$ with respect to $l$.
If $l_\infty=-1$, then $N_l=N^{\rho_f}$; 
if $l_\infty=2$, then $N_l/N^{\rho_f}=\{\mathrm{id},x\mapsto -x+4n\}
\subset\mathrm{Iso}(2\mathbb Z)$ for some $n\in\mathbb Z$.
Therefore if $\pi_r(\Sigma)\ne\{2\}$, $N_r=N^{\rho_f}$; otherwise, 
$N_r/N^{\rho_f}=\{\mathrm{id},x\mapsto -x+4n\}$.
In the latter case, the branched covering $\phi_{N_r}:S_{N_r}\to S$ is 
given by $\phi_{N_r}(z)=\phi(\sqrt z)=2\cos 2\pi\sqrt z$.
We have $\hat G=G^{\rho_f}$.

By Corollary \ref{co:34-1}, $$\mathrm{Cod}(f)\approx\{2\}\cup
\{[r]_{N^\rho}\in\mathrm{Rad}_{N^\rho}(f,\bar x)\,|\,\pi_r(\Sigma)
\ne\{2\}\}/\hat G,$$
where $2:\Sigma\to J$ is the coding map with image $\{2\}$.
Let us identify $\mathrm{Rad}_{N^\rho}(f,\bar x)
=(\phi^{-1}f^{-1}(\bar x))^2$ with $(\{\pm1\}\times \mathbb Z)^2$ 
by $\pm1/4+n\leftrightarrow (\pm1,n)$.
The action of $\mathrm{Iso}(2\mathbb Z)$ on 
$(\{\pm1\}\times\mathbb Z)^2$ 
is given by 
\begin{multline*}
(x\mapsto ax+2n)\cdot((\epsilon_1,n_1),(\epsilon_2.n_2))\\
=((\epsilon_1,an_1+(2-\epsilon_1)n),
(\epsilon_2,an_2+(2-\epsilon_2)n) .
\end{multline*}
It is easily seen that $\pi_r(\Sigma)=\{2\}$ if and only if
$[r]_{N^\rho}=(b_1,b_2)$, $b_1,b_2\in\{(1,n),(-1,3n)\}$ for some $n\in
\mathbb Z$.
Thus
\begin{multline*}
\mathrm{Cod}(f)\approx\bigl\{[(1,0),(1,0)]\bigr\}
\cup
\big\{[(\epsilon_1,0),(\epsilon_2,n)]\,:\,\epsilon_i\in\{\pm1\},
n\in\mathbb N\}\\
\cup\{[(-1,1),(\epsilon,m)]\,:\,\epsilon\in\{\pm1\},m\in\mathbb Z\big\}.
\end{multline*}

In general, consider $f_d(z)=2T_d(z/2)$, where 
$T_d$ is the Chebyshev polynomial of degree $d$ (i.e. $f_d\circ\phi(z)
=\phi(\pm dz+2n)$).
Then $f_2(z)=z^2-2$.
For the basepoint $\bar x=0$, $\phi^{-1}f_d^{-1}(\bar x)=\{(\pm1/2+2n)/d
\,|\,n\in\mathbb Z\}$.
Similarly to the above, 
$$\mathrm{Cod}(f_d)\approx\{2\}\,(\cup\{-2\})\cup
\{[r]_{N^\rho}\in\mathrm{Rad}_{N^\rho}(f_d,\bar x)\,|\,\pi_r(\Sigma)
\ne\{2\},\{-2\}\}/\hat G,$$ 
where the term $\{-2\}$ appears if $d$ is odd.
Letting $(\pm1/2+2n)/d$ correspond to $(\pm1,n)\in\{\pm1\}
\times\mathbb Z$, the action of $\hat G=\mathrm{Iso}(2\mathbb Z)$ 
on $\mathrm{Rad}_{N^\rho}(f,\bar x)\approx(\{\pm1\}\times\mathbb Z)^d$ 
is given by 
$$(x\mapsto ax+2n)\cdot((\epsilon_i,n_i))=((\epsilon_i,an_i
+(d-\epsilon_i)n)).$$

\item
$f(z)=-(z-1)^2/4z$.
The group $G^{\rho_f}$ is equal to $$\mathrm{Iso}^+(2\Gamma)=
\{z\mapsto a z+2b\,|\,a=\pm1,\pm i,\ b\in\Gamma\}.$$
Similarly to the above, if $\pi_r(\Sigma)\ne \{\infty\}$,  
$N_r=N^{\rho_f}$; otherwise $N_r/N^{\rho_f}
=\langle z\mapsto iz+2b\rangle$ for some $b\in\Gamma'$.
Hence $\hat G=G^{\rho_f}$.
We have $$\mathrm{Cod}(f)\approx\{\infty\}\cup
\{[r]_{N^\rho}\in\mathrm{Rad}_{\hat N}(f,\bar x)\,|\,\pi_r(\Sigma)
\ne\{\infty\}\}/\hat G,$$
where $\infty$ denotes the coding map with image $\{\infty\}$.
Letting $\alpha/2+(1+i)\beta\ (\alpha\in A=\{\pm1,\pm i\},
\beta\in\Gamma)$ correspond to $(\alpha,\beta)\in A\times \Gamma$, 
we have the action of $\mathrm{Iso}^+(2\Gamma)$  on $(A\times 
\Gamma)^2$ by 
\begin{multline*}
(z\mapsto a z+2b)\cdot((\alpha_1,\beta_1),
(\alpha_2,\beta_2))\\
=((\alpha_1,a\beta_1+(2-(1-i)\alpha_1) b),
(\alpha_2,a\beta_2+(2-(1-i)\alpha_2)b)).
\end{multline*}
Thus 
\begin{multline*}
 \mathrm{Cod}(f)=\{[(1,0),(1,0)]\}\cup\{[(\alpha_1,0),(\alpha_2,\beta)]
\,:\,\alpha_i\in A,\beta\in\Gamma_+\}\\
\cup\{[(\alpha_1,1+i),(\alpha_2,\beta)]\,:\,\alpha_1\in\{-1,-i\},
\alpha_2\in A,\beta\in\Gamma\},
\end{multline*}
where $\Gamma_+=\{n+mi\in\Gamma\,|\,n>0,m\ge 0\}$.
\end{enumerate}
\end{exa}

Let $\mathcal A(f)=\{R:\hat{\mathbb C}\to\hat{\mathbb C}
\,|\,R\text{ is rational, }\deg R\ge1, R\circ f=f\circ R\}$.
Then $\mathcal A(f)$ forms a monoid under composition of maps.

\begin{pro}\label{pro:32}
The monoid  $\mathcal A(f)$ acts on $\mathrm{Cod}(f)=\{\pi_r\}$ 
by $\pi\mapsto R\circ\pi$.
Moreover, $\mathrm{mul}(R\circ \pi)=\deg R\cdot\mathrm{mul}(\pi)$ 
for $\pi\in\mathrm{Cod}(f)$.
\end{pro}

\begin{proof}
For $R\in\mathcal A(f)$, $R(J)=J=R^{-1}(J)$.
Indeed, if $x$ belongs to the Fatou set of $f$, then so does $R(x)$ 
(see for example \cite{Be91}, Theorem 4.2.9).
Thus $J_f=J_R$ if $\deg R\ge 2$.
If $\deg R=1$, we have $J\subset R(J)$, so $R^{-1}(J)\subset J$.
On the other hand, $R(J)\subset J\subset R^{-1}(J)$ since $R^{-1}\in
\mathcal A(f)$.

Since $R:J\to J$ is $(\deg R)$-to-one except for finite points, 
$\mathrm{mul}(R\circ \pi)
=\deg R\cdot\mathrm{mul}(\pi)$.  
For a radial $r$, we can assume that 
$r$ does not intersect $R^{-1}(P_f)$. 
It is easily seen that $R\circ\pi_r=\pi_{R\circ r}$.
Thus $R\circ\pi_r\in\mathrm{Cod}(f)$.
\end{proof}

\begin{pro}\label{pro:40}
Every  $R\in \mathcal A(f)$ of degree $D=\deg R$
bigger than one is a subhyperbolic rational map with $P_f=P_R$ 
and $\rho_f=\rho_R$.
Moreover, exceptional points of $f$ coincide with those of $R$.
\end{pro}

\begin{proof}
Recall that  $p$ is called an exceptional point of $f$ if 
 $\#\bigcup_{k=0}^\infty f^{-k}(p)<\infty$.
In fact, then $p$ is a critical point of degree $d=\deg f$, $f^2(p)
=p$, $\bigcup_{k=0}^\infty f^{-k}(p)=\{p,f(p)\}$, and $f(p)$ is also 
exceptional.

First we show that if $x=R^k(y)$ and $y$ is not an exceptional point of 
$f$, then there exist $m>0$ and $z\in f^{-m}(x)$
such that $\deg_y(R^k)$ is a divisor of $\deg_{z}(f^m)$. 
Take $m>0$ so that $\#f^{-m}(y)>2D^k-2$.
Since $\#C_{R^k}\le 2D^k-2$, there exists $w\in f^{-m}(y)$ such that 
$\deg_w(R^k)=1$.
Setting $z=R^k(w)$, 
we have $\deg_y(R^k)\deg_w(f^m)=\deg_w(R^kf^m)=\deg_w(f^mR^k)
=\deg_{z}(f^m)\deg_{w}(R^k)
=\deg_{z}(f^m)$.

Second we show that if $y$ is an exceptional point of $f$, then $y$ is an 
exceptional point of $R$.
Suppose $y$ is an exceptional point of $f$.
By commutativity, $R(y),R^2(y),\dots$ are exceptional points of $f$.
Since a rational map has at most two exceptional points, $R^3(y)=R(y)$ 
and $f^{-2}(R(y))=\{R(y)\}$.
Assume that $R(y)$ is not an exceptional point of $R$.
Then $X_k=R^{-2k}(R(y))-\{R(y)\}$ is nonempty for $k\in\mathbb N$.
From $R^{2k}\circ f^2(X_k)=R(y)$, we have $f^2(X_k)\subset X_k$.
Hence $X_k$ contains a periodic cycle of $f$.
This contradicts the fact that
 the Fatou set includes at most finite periodic cycles.

From the claims above, we deduce that $P_f=P_R$ and $\rho_f=\rho_R$. 
\end{proof}

\begin{co}
 If $f$ has exactly one critical point $c$ with $\#\{f^k(c)\,|\,
k\in\mathbb N\}=\infty$, then $\mathcal A(f)$ is generated by $f$.
\end{co}

\begin{proof}
Let $R\in\mathcal A(f)$.
We show that $R\circ f(c)=f^k(c)$ for some $k\in\mathbb N$.
To this end, we assume that 
$R\circ f(c)\notin\{f^k(c)\,|\,k\in\mathbb N\}$. 
Since $f(c)\in P$, $R\circ f(c)\in P$.
Thus $R\circ f(c)$ is eventually periodic under the iteration of $f$.
Therefore there exist $m,n\in\mathbb N$ such that 
$f^m\circ R\circ f(c)=f^{m+kn}\circ R\circ f(c)=R\circ 
f^{m+kn+1}(c)$ for every $k\in\mathbb N$.
This implies 
a contradiction that $\#R^{-1}(f^m\circ R\circ f(c))=\infty$. 

From $R\circ f(c)=f^k(c)$, we have $R(f^n(c))=f^{n+k-1}(c)$.
Hence $R=f^{k-1}$ by the identity theorem.
\end{proof}

\begin{co}
Let $\rho=\rho_f$ be the canonical ramification function, and 
$g:S_{N^\rho}\to S_{N^{\rho}}$ a contraction with $f\circ \phi_{N^\rho}
\circ g=\phi_{N^\rho}$.
Then $\mathcal A(f)$ is the set of rational maps $R:\hat{\mathbb C}
\to\hat{\mathbb C}$ such that there exists 
a branched covering $g':S_{N^\rho}\to S_{N^\rho}$ with 
$R\circ \phi_{N^\rho}\circ g'=\phi_{N^\rho}$ and 
$g\circ g'=g'\circ g\circ t$ for some deck transformation $t$.
\end{co}

\begin{proof}
If $R\in\mathcal A(f)$, then  there exists 
a branched covering $g':S_{N^\rho}\to S_{N^\rho}$ with 
$R\circ \phi_{N^\rho}\circ g'=\phi_{N^\rho}$ by Proposition \ref{pro:40}.
From $f\circ R\circ \phi_{N^\rho}\circ g'\circ g=\phi_{N^\rho}$ 
and $f\circ R\circ \phi_{N^\rho}\circ g\circ g'
=R\circ f\circ \phi_{N^\rho}\circ g\circ g'=\phi_{N^\rho}$, it follows 
 that $g\circ g'=g'\circ g\circ t$ for some $t$ by Proposition 
\ref{pro:1}.  

Conversely, suppose there exists 
a branched covering $g':S_{N^\rho}\to S_{N^\rho}$ with 
$R\circ \phi_{N^\rho}\circ g'=\phi_{N^\rho}$ and 
$g\circ g'=g'\circ g\circ t$ for some deck transformation $t$.
Then $R\circ f\circ \phi_{N^\rho}\circ g\circ g'=\phi_{N^\rho}$
and $f\circ R\circ \phi_{N^\rho}\circ g'\circ g\circ t=\phi_{N^\rho}$.
Hence
\begin{eqnarray*}{}
R\circ f(x)&=&R\circ f\circ\phi_{N^\rho}\circ g\circ g'
((\phi_{N^\rho}\circ g\circ g')^{-1}(x))\\
&=&\phi_{N^\rho}((\phi_{N^\rho}\circ g\circ g')^{-1}(x))
\ =\ \phi_{N^\rho}((\phi_{N^\rho}\circ g'\circ g'\circ t)^{-1}(x))\\
&=&f\circ R\circ\phi_{N^\rho}\circ g\circ g'\circ t
((\phi_{N^\rho}\circ g'\circ g\circ t)^{-1}(x))
\ =\ f\circ R(x)
\end{eqnarray*}
for any $x\in S$.
\end{proof}

\begin{defi}
 A coding map $\pi\in \mathrm{Cod}(f)$ is said to be {\em prime} 
if there are no $R\in\mathcal A(f)$ and no $\pi'\in\mathrm{Cod}(f)$ 
such that $\pi=R\circ \pi'$ and $\deg R\ge 2$.
A radial $r$ is said to be prime if $\pi_r$ is prime.

We write $$
\begin{array}{l}
\mathrm{Cod}'(f)=\{\pi\in\mathrm{Cod}(f)\,|\,\text{$\pi$ is prime}\}
\approx\mathrm{Cod}(f)/\mathcal A(f)\\
\mathcal M=\{\mathrm{mul}(\pi)\,|\,\pi\in\mathrm{Cod}(f)\},\ 
\mathcal M'=\{\mathrm{mul}(\pi)\,|\,\pi\in\mathrm{Cod}'(f)\}
\end{array}$$
\end{defi}

\begin{co}
{\rm(1)} $0\in\mathcal M'$,
{\rm (2)} $\mathcal M'-\{0\}\ne\emptyset$,
{\rm(3)} $\mathcal M=\{\deg R\cdot n\,|\,R\in\mathcal A(f),n\in\mathcal 
M'\}$. 
In particular, $\mathcal M\supset
\{d^k\cdot n\,|\,n\in\mathcal M',k=0,1,2,\dots\}$ since $f^k\in\mathcal 
A(f)$.  
\end{co}

\begin{exa}
 We calculate $\mathcal M'$, $\mathcal M$ and $\mathrm{Cod}'(f)$ 
for a couple of examples in Section \ref{sec:3}.
\begin{enumerate}
 \item  $f(z)=z^d$.
The monoid 
$\mathcal A(f)=\{e^{2\pi im/d}z^k\,|\,k\in\mathbb Z,m=0,1,\dots,d-1
\}$ is identified with 
$\{(m,k)\}$ (the product is given by $(m,k)(m',k')=(m+km',kk')$).
The action of $\mathcal A(f)$ on $\mathrm{Cod}(f)=\mathbb Z^d/\mathbb 
Z$ is given by 
$$(m,k)\cdot[n_1,n_2,\dots,n_d]=[kn_1+m,kn_2+m,\dots,kn_d+m].$$
Hence $$\mathrm{Cod}(f)/\mathcal A(f)\approx\{[0,n_2,\dots,n_d]\,|\,
n_2,\dots,n_d\text{ are mutually prime}\}.$$
Since $1\in\mathcal M'$, we have $\mathcal M=\{0\}\cup\mathbb N$.

In the case $d=2$, 
$\mathrm{Cod}(f)$ is identified with $\mathbb Z$ by $n_2-n_1
\leftrightarrow [n_1,n_2]$, and $\mathrm{mul}([n_1,n_2])=n_2-n_1$. 
We can see that $\mathrm{Cod}'(f)$ is identified with 
$\{[0,0],[0,1],[0,-1]\}$.
Thus $\mathcal M'=\{0,1\}$.

In the case $d=3$, $\mathrm{Cod}'(f)\approx
\{[0,n_2,n_3]\,|\,\text{$n_2$, 
$n_3$ are mutually prime}\}$ is infinite, 
and $\mathcal M'=\{0,1\}$.

In the case $d=4$, $\mathcal M'\supset\{0,1,2\}$ 
(see \cite{LaWa96}, Example 3.1).

\item $f(z)=z^2-2$.
Then $\mathrm{Cod}(f)\approx(\{+,-\}\times \mathbb Z)^2/
\mathrm{Iso}(\mathbb Z)$ and $\mathcal A(f)=\{f_d\,|\,d=1,2,\dots\}
\approx\mathbb N$, where $f_d(z)=2T_d(z/2)$.
The action of $\mathcal A(f)$ on $\mathrm{Cod}(f)$ is given by 
$n\cdot[(\epsilon_1,n_1),(\epsilon_2,n_2)]=[(\epsilon_1,nn_1),
(\epsilon_2,nn_2)]$.
Hence 
\begin{multline*}
\mathrm{Cod}'(f)\approx\big\{[(+,0),(+,0)],[(-,1),(-,1)],
[(-,1),(+,0)]\big\}\\
\cup\big\{[(\epsilon_1,0),(\epsilon_2,1)]\,:\,\epsilon_i\in\{+,-\}
\big\}
\cup\big\{[(-,1),(-,\pm3^n+1)]\,:\,n\in\mathbb N\cup\{0\}\big\}.
\end{multline*}
Thus $\mathcal M'=\{0,1,2\cdot3^n\,|\,n=0,1,\dots\}$ 
and $\mathcal M=\{0\}\cup
\mathbb N$.

\item $f(z)=-(z-1)^2/4z$.
The monoid $\mathcal A(f)=\{R_b\,|\,b\in\Gamma\}\approx\Gamma$, 
where $R_b$ is the rational map with 
$R_b\circ\phi(z)=\phi(b z)$.
The action of $\mathcal A(f)$ on $\mathrm{Cod}(f)$ is given by 
$b\cdot[(\alpha_1,\beta_1),(\alpha_2,\beta_2)]=[(\alpha_1,b\beta_1),
(\alpha_2,b\beta_2)]$.
Hence 
\begin{multline*}
\mathrm{Cod}'(f)\approx\big\{[(1,0),(1,0)]\big\}
\cup\big\{[(-1,1),(-i,1+i)],[(-i,1),(-1,1+i)]\big\}\\
\cup\big\{[(\alpha_1,0),(\alpha_2,1)]\,:\,(\alpha_1,\alpha_2)\in A^2
-\{(-1,-i),(-i,-1)\}\big\}\\
\cup\big\{[(\alpha_1,1),(\alpha_2,0)]\,:\,(\alpha_1,\alpha_2)
\in\{-1,-i\}\times A-\{(-1,-i),(-i,-1)\}\big\}\\
\cup\big\{[(\alpha,1),(\alpha,1)]\,:\,\alpha\in\{-1,-i\}\big\}\\
\cup\big\{[(-1,1),(-1,(1-2i)^n\alpha+1)]\,:\,\alpha\in A,
n\in\mathbb N\cup\{0\}\big\}\\
\cup\big\{[(-i,1),(-i,(1+2i)^n\alpha+1)]\,:\,\alpha\in A,
n\in\mathbb N\cup\{0\}\big\}.
\end{multline*}
Thus $\mathcal M'=\{0,1,4\cdot5^n\,|\,n=0,1,\dots\}$ 
and $\mathcal M=\{|z|^2\,:\,z\in\Gamma\}$.

\item $f(z)=z^2-3$.
Then $\mathcal A(f)$ is generated by $f$.
In this case, the complete solution has not been obtained.
We will show that $1,2\in \mathcal M'$ in Example \ref{exa:35} 
(we conjecture $\mathcal M'=\{1,2\}$).
Consequently,  $0,2^k\in \mathcal M,k\in \mathbb N$.
\end{enumerate}
 \end{exa}

We pose the following problem:

\medskip

\noindent{\bf Problem.}
For a given $f$, determine $\mathcal M'$ and $\mathcal M$.

\medskip

It is easily seen that $\mathcal M=\{0\}\cup\mathbb N$ 
if $f$ is postcritically finite polynomial map.
In general case, however, it is unknown whether $1\in\mathcal M$ or not.

\section{Equivalence relations on the word space}\label{sec:5}

 Let $r=(l_i)$ and $r'=(l_i')$ be radials with basepoints $\bar x$ and 
$\bar x'$ respectively, and $N\subset G$ a subgroup.
For $r'$, the notation $x_w'$, $l'_w$, $l'_\omega$ and $F'_x(\cdot)$ 
are defined in a trivial way.
 Let $M$ be a real number bigger than $\sup_{\omega\in\Sigma}|l_\omega| 
+\sup_{\omega\in\Sigma}|l'_\omega|$, and let 
$$Y'=\{\gamma\in Y\,:\,|\gamma|<M\},\quad Y_N'=Y'/\sim_N,$$
where $\gamma\sim_N \gamma'$ if $[\gamma{\gamma'}^{-1}]\in  N$.
It is clear that $Y_N'$ is a finite set if $N^{\rho_f}\subset N$.
The equivalence class of $\gamma$ is denoted by $[\gamma]_{N}.$

\begin{pro}\label{pro:16}
For $\omega=\omega_1\omega_2\cdots,\omega'=\omega_1\omega_2'\cdots
\in\Sigma$, $\pi_r(\omega)=\pi_{r'}(\omega')$ 
if and only if there exist curves $\gamma_0,\gamma_1,\ldots \in Y'$ 
such that 
\begin{equation}
[l_{\omega_k}F_{x_{\omega_k}}(\gamma_k)
{l'_{\omega_k'}}^{-1}\gamma_{k-1}^{-1}]\text{ is trivial}\label{eq:4}
\end{equation}
 for $k=1,2,\dots$.
Moreover, {\rm (\ref{eq:4})} can be replaced with 
\begin{equation}
[l_{\omega_k}F_{x_{\omega_k}}(\gamma_k)
{l'_{\omega_k'}}^{-1}\gamma_{k-1}^{-1}]\in{ N_r}.\label{eq:5}
\end{equation}
\end{pro}

\begin{proof}
 Suppose $p=\pi_r(\omega)=\pi_{r'}(\omega')$.
We take $\gamma_k=l_{\sigma^k\omega}{l'_{\sigma^k\omega'}}^{-1}$.
Then $F_{x_{\omega_k}}(\gamma_k)=F_{x_{\omega_k}}(l_{\sigma^k(\omega)})
{{F'}_{x'_{\sigma^k\omega'}({l'}_{\sigma^k\omega})}}^{-1}$.
Thus if $p\notin P$, then $\gamma_0,\gamma_1,\dots$ satisfy the condition 
above.
If $p\in P$, then modify each $\gamma_k$ in a small neighborhood of
 $f^k(p)$ so that (\ref{eq:4}) holds.

Conversely, suppose there exist curves $\gamma_0,\gamma_1,\dots$ 
satisfying the condition (\ref{eq:5}).
We write $(\omega)_k=\omega_1\omega_2\cdots\omega_k$.
Taking the product of 
$$\begin{array}{l}
\ [l_{(\omega)_k}F_{x_{(\omega)_k}}(\gamma_k){l'_{(\omega')_k}}^{-1}
l'_{(\omega')_{k-1}}F_{x_{(\omega)_{k-1}}}(\gamma_{k-1})^{-1}
l_{(\omega)_{k-1}}^{-1}] \\
\qquad\qquad=[l_{(\omega)_{k-1}}F_{x_{(\omega)_{k-1}}}
(l_{\omega_k}F_{x_{\omega_k}}(\gamma_k)
{l'_{\omega_k'}}^{-1}\gamma_{k-1}^{-1})
l_{(\omega)_{k-1}}^{-1}]\in N_r
\end{array}$$
from $k=n$ to $1$,  we have 
$[l_{(\omega)_n}F_{x_{(\omega)_n}}(\gamma_n){l'_{(\omega)_n}}^{-1}
\gamma_0]\in{ N_r}.$
Thus the distance between $x_{(\omega)_n}$ and $x'_{(\omega')_n}$ 
is less than $|F_{x_{(\omega)_n}}(\gamma_n)|\le c^{-n}M$.
\end{proof}

\begin{thm}\label{thm:28}
There exists a weighted directed graph $(V,E,\alpha)$:
\begin{itemize}
\item The vertex set $V$ is finite.
\item The edge set $E$ is finite. 
Each edge $e\in E$ has its initial vertex $e^-\in V$ and its terminal 
vertex $e^+\in V$.
(We do not assume that $e_0^-=e_1^-$ and $e_0^+=e_1^+$ imply $e_0=e_1$.)
\item The weight function $\alpha:E\to\{1,2,\dots,d\}^2$.  
\end{itemize}
such that for $\omega=\omega_1\omega_2\cdots,\omega'=\omega_1'\omega_2'
\cdots\in\Sigma$, 
$$\pi_r(\omega)=\pi_{r'}(\omega') \iff
\begin{array}{l}
\text{there exist $e_1,e_2,\cdots\in E$ such that } \\
\text{$e_i^+=e_{i+1}^-$ and $\alpha(e_i)=(\omega_i,\omega_i') $.}
\end{array}$$
\end{thm}

\begin{proof}
Let $V$ be the maximal subset of $Y_{N_r}'$ such that for any 
$[\gamma]_{N_r}\in Y_{N_r}'$ 
there exists $[\gamma']_{N_r}\in V$ with $[\gamma]_{ N_r}
=[l_iF_{x_i}(\gamma'){l_j'}^{-1}]_{N_r}$ 
for some $i,j\in\{1,2,\dots,d\}$.
Set 
$$E=\{([\gamma]_{N_r},[\gamma']_{N_r},i,j)\in V^2\times\{1,2,\dots,d\}^2 
\,|\,[\gamma]_{N_r}=[l_iF_{x_i}(\gamma'){l_j'}^{-1}]_{N_r}\}.$$
For $e=([\gamma]_{ N_r},[\gamma']_{ N_r},i,j)\in E$, 
define $e^-=[\gamma]_{N_r}$, 
$e^+=[\gamma']_{N_r}$ and $\alpha(e)=(i,j)$.

Suppose there exist $e_1,e_2,\cdots\in E$ such that 
$e_i^+=e_{i+1}^-$ and $\alpha(e_i)=(\omega_i,\omega_i')$.
Since there exist curves $\gamma_k,k=0,1,\dots$ such that 
$[\gamma_k]_{N_r}
=e_{k+1}^-$, we have $\pi_r(\omega)=\pi_{r'}(\omega')$ by proposition 
\ref{pro:16}.
Conversely, suppose $\pi_r(\omega)=\pi_{r'}(\omega')$.
By Proposition \ref{pro:16}, there exists $[\gamma_0]_{N_r},
[\gamma_1]_{N_r},\dots\in Y'$ such that $[l_{\omega_k}
F_{x_{\omega_k}}(\gamma_k)
{l'_{\omega_k'}}^{-1}]_{N_r}=[\gamma_{k-1}]_{N_r}$ 
for $k=1,2,\dots$.
Then $[\gamma_k]_{N_r}\in V$, and so $e_k=([\gamma_{k-1}]_{N_r},
[\gamma_k]_{N_r},\omega_k,\omega_k'),k=1,2,\dots$ 
satisfy the condition.
\end{proof}

\begin{co}\label{co:34}
 Let $V$ be the vertex set constructed in Theorem \ref{thm:28}.
If $\pi(\Sigma)=J$, then $$\#\pi^{-1}(x)\le 
\max_w\#\{w'\,|\,x_w=x_{w'}\} \cdot \max_{p\in C\cap J,k\ge 1}
\deg_p f^k\cdot\sqrt{2\#V}$$
for any $x\in J$.
\end{co}

\begin{proof}
 An immediate consequence of  (\ref{eq:l}).
\end{proof}

Is is easily seen that if $\pi(\Sigma)=J$, then 
the Julia set $J$ is topologically identified 
with the quotient space $\Sigma/\sim_r$, where the equivalence relation 
$\sim_r$ is defined by $\omega\sim_r\omega'\iff \pi_r(\omega)
=\pi_r(\omega')$.
Considering $r=r'$, we obtain from Theorem \ref{thm:28} 
 an algorithm to calculate $\sim_r$ 
provided $Y'$ is determined. 

\begin{exa}\label{exa:35}
 Let $f(z)=z^2-3$.
In order to obtain $\sim_r$ for the radials $r=r_1,r_2,r_3$ in Section 
\ref{sec:3.1}, we use Theorem \ref{thm:28}.
Let us set generators $[B_1],[B_2],\dots$ of $G$ as shown in 
Figure \ref{fig:fig5}.
Then $[B_k^2],[(B_1B_2)^4]\in\hat N,k=1,2,\dots$
Take a simply connected open domains $U_1\subset U_2$ 
such that $B_1\subset U_1,B_2\subset U_2,
[-3,3]\subset U_1,[-3,6]\subset U_2,
[6,\infty)\subset\mathbb C-U_1,[33,\infty)\subset \mathbb C-U_2$, 
and $f^{-1}(\overline {U_i})\subset U_i$.
For the radial $r=r_j$ above, we take $r$ so that the image of $r$ 
is included in $U_i$ ($i=1$ if $j=1,2$, $i=2$ if $j=3$) 
without loss of generality.
For any $k$, there exists $n$ such that 
$F_x(B_k)\subset U_2$ for each $x\in f^{-n}(\bar x)$.
Therefore, for $r=r_1,r_2$, 
the set $V$ in Theorem \ref{thm:28} is included 
$\{1_{\hat N},[B_1]_{\hat N}\}$, 
where $1_{\hat N}=\hat N$ is the unit element.
If $r=r_3$, $V$ is included in 
\begin{multline*}
\{1_{\hat N},[B_1]_{\hat N},[B_2]_{\hat N},[B_1B_2]_{\hat N},
[B_2B_1]_{\hat N},[B_1B_2B_1]_{\hat N},[B_2B_1B_2]_{\hat N},\\
[(B_1B_2)^2]_{\hat N}[(B_2B_1)^2]_{\hat N}\},
\end{multline*}
the subgroup generated by $[B_1]_{\hat N}$ and $[B_2]_{\hat N}$.
\begin{enumerate}
 \item $r=r_1$.
We have 
$$\begin{array}{lcll}
e&\leftarrow&(11)\ e,&(22)\ e\\
B_1&\leftarrow& (12)\ e,&(21)\ e
\end{array}$$
 where $e$ denotes the trivial loop.
(For example, ``$B_1\leftarrow(12)\ e,\ (21)\ e$'' indicates that 
$l_1F_{x_1}(B_1)l_2^{-1}$ and $l_2F_{x_2}(B_1)l_1^{-1}$ are  
 trivial.)
Thus $$V=\{1_{\hat N}\},\ E=\{(1_{\hat N},1_{\hat N},1,1),
(1_{\hat N},1_{\hat N},2,2)\}.$$
It follows from this that the equivalence relation $\sim_r$ is trivial,
that is, $\omega\sim_r\omega'$ if and only if $\omega=\omega'$.
Consequently, $\pi_r$ is bijective.
\item $r=r_2$.
We have 
$$
\begin{array}{lcll}
e&\leftarrow&(11)\ e,&(22)\ e\\
B_1&\leftarrow& (12)\ B_1^{-1},&(21)\ B_1\\
\end{array}$$
(``$B_1\leftarrow(11)\ B_1^{-1},\ (22)\ B_1$'' indicates that 
$l_1F_{x_1}(B_1)l_2^{-1}$ is homotopic to $B_1^{-1}$ and 
$l_2F_{x_2}(B_1)l_1^{-1}$ is homotopic to $B_1$).
Thus $$\begin{array}{l}V=\{1_{\hat N},[B_1]_{\hat N}\},\\
 E=\{(1_{\hat N},1_{\hat N},1,1),(1_{\hat N},1_{\hat N},2,2),
([B_1]_{\hat N},[B_1]_{\hat N},1,2),
([B_1]_{\hat N},[B_1]_{\hat N},2,1)\}.
\end{array}$$
It follows from this that 
$$\omega=\omega_1\omega_2\cdots\sim_r \omega'=\omega_1'\omega_2'\cdots
 \iff \omega=\omega'\text{ or }\omega_k\ne\omega_k',k=1,2,\dots$$
(for example, $111\cdots\sim_r222\cdots$ and 
 $1212\cdots\sim_r2121\cdots$). 
Consequently, $\pi_r$ is exactly two-to-one.
\item $r=r_3$.
We have 
$$
\begin{array}{lcll}
e&\leftarrow&(11)\ e,&(22)\ e\\
B_1&\leftarrow& (12)\ B_2,&(21)\ B_2^{-1}\\
B_2&\leftarrow& (11)\ e,&(22)\ B_2^{-1}B_1B_2\\
B_1B_2&\leftarrow& (12)\ B_1B_2,&(21)\ B_2^{-1}\\
B_2B_1&\leftarrow& (12)\ B_2,&(21)\ B_2^{-1}B_1\\
B_2B_1B_2&\leftarrow& (12)\ B_1B_2,&(21)\ B_2^{-1}B_1
\end{array}$$
etc.
Thus 
$$\begin{array}{ccl}
 V&=&\{1_{\hat N}, [B_2]_{\hat N},[B_1B_2]_{\hat N},[B_2B_1]_{\hat N},
[B_2B_1B_2]_{\hat N}\},\\
 E&=&\{(1_{\hat N},1_{\hat N},1,1),(1_{\hat N},1_{\hat N},2,2),\\
&&\qquad(1_{\hat N},[B_2]_{\hat N},1,1),
([B_2]_{\hat N},[B_2B_1]_{\hat N},1,2),\\
&&\qquad ([B_2]_{\hat N},[B_1B_2]_{\hat N},2,1),
([B_2B_1]_{\hat N},[B_2B_1]_{\hat N},2,1),\\
&&\qquad([B_1B_2]_{\hat N},[B_1B_2]_{\hat N},1,2),
 ([B_2B_1]_{\hat N},[B_2B_1B_2]_{\hat N},2,1),\\
&&\hspace{1.5cm}([B_1B_2]_{\hat N},[B_2B_1B_2]_{\hat N},1,2),
([B_2B_1B_2]_{\hat N},[B_2]_{\hat N},2,2)\}.
\end{array}$$
By Corollary \ref{co:34}, $\pi_r$ is at most three-to-one.
The graph $(V,E,\alpha)$ is diagrammatically shown as
$$\begin{array}{crclc}
{\scriptstyle 11,22}\ \text{\scalebox{1.5}{\rotatebox[origin=c]{270}{$\circlearrowleft$}}\ }
1&\xrightarrow{\ 11\ } \quad B_2\ &\xrightarrow{\ 12\ }&B_2B_1\
 \text{\scalebox{1.5}{\rotatebox[origin=c]{90}{$\circlearrowleft$}}}
\ {\scriptstyle 21}\\
&{\scriptstyle 21}\downarrow\ \ \ & \ \nwarrow{\scriptstyle 22}   &
\ \ \downarrow{\scriptstyle 21}      \\
&{\scriptstyle 12}
\ \text{\scalebox{1.5}{\rotatebox[origin=c]{270}{$\circlearrowleft$}}\ }
B_1B_2& \xrightarrow[\ 12\ ]{}& B_2B_1B_2
\end{array}$$
We can see that if $\omega\in\Sigma$ contains a word $12121$ 
 infinitely many times, then $\pi_r^{-1}\pi_r(\omega)=\{\omega\}$.
Consequently, the multiplicity of $\pi_r$ is equal to one.
\end{enumerate}
\begin{figure}[hbtp]
\centering
\includegraphics[height=5cm]{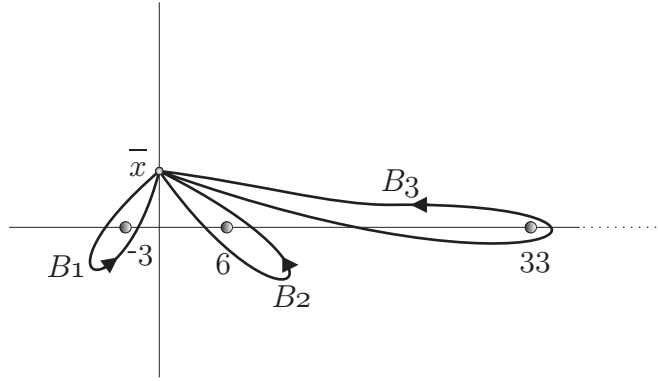}
\caption{Generators of $G$.}
\label{fig:fig5}
\end{figure}

\end{exa}

\section{Non-subhyperbolic case}\label{sec:nonhyp}
A couple of our results are true for non-subhyperbolic rational maps.
In fact, if $f$ is geometrically finite (i.e. $J\cap P$ is finite),
then almost all of our results are applicable, but 
it is possible that $J\not\subset S$.
Thus the Julia tile $K$ might be noncompact.
The details are left to the reader.

For general rational maps, we need some restriction. 
Let $f:\hat{\mathbb C}\to\hat{\mathbb C}$ be a rational map.
If $S'=\hat{\mathbb C}-P$ has a connected component $U$ such that 
$f^{-1}(U)\subset U$, then we have 
a radial $r$ in $U$,
but $x_\omega$ does not converge to a point in general.
We define
 $$\Pi(\omega)=\bigcap_{k=1}^\infty \overline{
\bigcup_{\omega'\in\Sigma}
F_{x_{\omega_1\omega_2\cdots\omega_{k}}}(l_{\omega'})}$$
for $\omega=\omega_1\omega_2\cdots\in\Sigma$.
Then we have

\begin{pro}
 {\rm(1)} $f(\Pi(\omega))=\Pi({\sigma\omega})$,
{\rm (2)} $\Pi(\omega)\subset J$,  
{\rm (3)} $\omega\mapsto \Pi(\omega)$ is upper semicontinuous, that is,  
if a sequence $\omega^1,\omega^2,\dots$ in $\Sigma$ converges to 
$\omega$, then $\bigcap_{m=1}^\infty\overline{\bigcup_{k=m}^\infty
\Pi({\omega^k})}\subset \Pi(\omega)$ 
 and {\rm (4)} $\Pi(\omega)$ is connected.
\end{pro}

\begin{proof}
(1) is immediately deduced from the definition.

(2): Suppose a point $x\in\Pi(\omega)$ belongs to the Fatou set.
It is easily seen that $f^n(x)$ converges to neither a 
(super)attracting cycle nor a parabolic cycle.
Thus $f^n(x)\in \Pi(\sigma^n\omega)$ for some $n$ is contained in 
either a Siegel disc or a Herman ring. 
This contradicts the fact that  $f:f^{-1}(U)\to U$ is expanding 
in the Poincar\'e metric on $U$.  

(3): Set
 $\omega=\omega_1\omega_2\cdots$ and $\omega^k=\omega_1^k\omega_2^k
\cdots$.
For any $n$, there exists 
$m_0$ such that $\omega^m\in \Sigma(\omega_1\omega_2\cdots\omega_n)$ 
whenever $m\ge m_0$.
Since $$\Pi(\omega^m)\subset
\overline{\bigcup_{\omega'\in\Sigma}
F_{x_{\omega_1^m\omega_2^m\cdots\omega_{n}^m}}(l_{\omega'})}=
\overline{\bigcup_{\omega'\in\Sigma}
F_{x_{\omega_1\omega_2\cdots\omega_{n}}}(l_{\omega'})}$$ for $m\ge m_0$,
we have $\bigcap_{m=1}^\infty\overline{\bigcup_{k=m}^\infty
\Pi({\omega^k})}\subset 
\overline{\bigcup_{\omega'\in\Sigma}
F_{x_{\omega_1\omega_2\cdots\omega_{n}}}(l_{\omega'})}$ for every $n$. 

(4): $\Pi(\omega)$ is connected since 
$\overline{\bigcup_{\omega'\in\Sigma}
F_{x_{\omega_1\omega_2\cdots\omega_{k}}}(l_{\omega'})}$ is connected.
\end{proof} 

From (2) and (4) immediately, 

\begin{co}
If $J$ is totally disconnected, then $\Pi(\omega)$ is a singleton for 
every $\omega\in\Sigma$.
\end{co}

\begin{rem}
In general, by Przytycki's result \cite{Prz86},  
$\Pi(\omega)$ is a singleton for every $\omega\in\Sigma$ 
except for $\omega$ in a ``thin'' set.
\end{rem}

\begin{thm}
The multiplicity of $\Pi$ is well-defined, that is, there exists 
$n_r\in\mathbb N\cup \{\infty\}$ such that 
$\#\{\omega\in\Sigma\,|\,x\in\Pi(\omega)\}=n_r$ 
for $\mu$-almost all $x\in J$.
\end{thm}

\begin{proof}
It is sufficient to show that the function 
$h(x)=\#\{\omega\in\Sigma\,|\,x\in\Pi(\omega)\}$ is Borel. 
Let $A_{k,\epsilon}$ be the set of $x\in J$ 
such that for some $\omega^i\in K,i=1,2,\dots,k$, 
$x\in \Pi(\omega^i),1\le i\le k$ and
the distance between $\omega^i$ and $\omega^j$ is  equal to or bigger 
than $\epsilon$ for $0\le i\ne j\le k$, where we consider an  
arbitrary compatible distance function on $\Sigma$. 
Then $A_{k,\epsilon}$ is closed since $\omega\mapsto\Pi(\omega)$ is 
upper semicontinuous.
Thus $\{h(x)\ge k\}$ is Borel.
\end{proof}

\end{document}